\def\DATE{{November 14, 2008}}
\def\DATE{{\relax}}
\magnification=1100
\baselineskip=14truept
\voffset=.75in  
\hoffset=1truein
\hsize=4.5truein
\vsize=7.75truein
\parindent=.166666in
\pretolerance=500 \tolerance=1000 \brokenpenalty=5000

\footline={\hfill{\rm\the\pageno}\hfill\llap{\sevenrm\DATE}}

\def\note#1{%
  \hfuzz=50pt%
  \vadjust{%
    \setbox1=\vtop{%
      \hsize 3cm\parindent=0pt\eightpoints\baselineskip=9pt%
      \rightskip=4mm plus 4mm\raggedright#1%
      }%
    \hbox{\kern-4cm\smash{\box1}\hfil\par}%
    }%
  \hfuzz=0pt
  }
\def\note#1{\relax}

\newcount\equanumber
\equanumber=0
\newcount\sectionnumber
\sectionnumber=0
\newcount\subsectionnumber
\subsectionnumber=0
\newcount\snumber  
\snumber=0

\def\section#1{%
  \subsectionnumber=0%
  \snumber=0%
  \equanumber=0%
  \advance\sectionnumber by 1%
  \noindent{\bf \the\sectionnumber .~#1~}%
}%
\def\subsection#1{%
  \advance\subsectionnumber by 1%
  \snumber=0%
  \equanumber=0%
  \noindent{\bf \the\sectionnumber .\the\subsectionnumber .~#1.~}%
}%
\def\prevs{\the\sectionnumber .\the\subsectionnumber .\the\snumber }

\long\def\Corollary#1{%
  \global\advance\snumber by 1%
  \bigskip
  \noindent{\bf Corollary~\prevs .}%
  \quad{\it#1}%
}%
\long\def\Lemma#1{%
  \global\advance\snumber by 1%
  \bigskip
  \noindent{\bf Lemma~\prevs .}%
  \quad{\it#1}%
}%
\def\Proof{\noindent{\bf Proof.~}}
\long\def\Proposition#1{%
  \advance\snumber by 1%
  \bigskip
  \noindent{\bf Proposition~\prevs .}%
  \quad{\it#1}%
}%
\long\def\Theorem#1{%
  \advance\snumber by 1%
  \bigskip
  \noindent{\bf Theorem~\prevs .}%
  \quad{\it#1}%
}%
\long\def\Statement#1{%
  \advance\snumber by 1%
  \bigskip
  \noindent{\bf Statement~\prevs .}%
  \quad{\it#1}%
}%
\def\ifundefined#1{\expandafter\ifx\csname#1\endcsname\relax}
\def\labeldef#1{\global\expandafter\edef\csname#1\endcsname{\prevs}}
\def\labelref#1{\expandafter\csname#1\endcsname}
\def\label#1{\ifundefined{#1}\labeldef{#1}\note{$<$#1$>$}\else\labelref{#1}\fi}

\def\preveq{(\the\sectionnumber .\the\subsectionnumber .\the\equanumber)}
\def\neq{\global\advance\equanumber by 1\eqno{\preveq}}

\def\ifundefined#1{\expandafter\ifx\csname#1\endcsname\relax}

\def\equadef#1{\global\advance\equanumber by 1%
  \global\expandafter\edef\csname#1\endcsname{\preveq}%
  \preveq}

\def\equaref#1{\expandafter\csname#1\endcsname}
 
\def\equa#1{%
  \ifundefined{#1}%
    \equadef{#1}%
  \else\equaref{#1}\fi}

\font\eightrm=cmr8%
\font\sixrm=cmr6%

\font\eightsl=cmsl8%

\font\eightbf=cmb8%

\font\eighti=cmmi8%
\font\sixi=cmmi6%

\font\eightsy=cmsy8%
\font\sixsy=cmsy6%

\font\eightex=cmex8%
\font\sixex=cmex6%
\font\fiveex=cmex5%

\font\eightit=cmti8%

\font\eighttt=cmtt8%

\font\tenbb=msbm10%
\font\eightbb=msbm8%
\font\sevenbb=msbm7%
\font\sixbb=msbm6%
\font\fivebb=msbm5%
\newfam\bbfam  \textfont\bbfam=\tenbb  \scriptfont\bbfam=\sevenbb  \scriptscriptfont\bbfam=\fivebb%

\font\tenbbm=bbm10

\font\tencmssi=cmssi10%
\font\sevencmssi=cmssi7%
\font\fivecmssi=cmssi5%
\newfam\ssfam  \textfont\ssfam=\tencmssi  \scriptfont\ssfam=\sevencmssi  \scriptscriptfont\ssfam=\fivecmssi%

\font\tenfrak=cmfrak10%
\font\eightfrak=cmfrak8%
\font\sevenfrak=cmfrak7%
\font\sixfrak=cmfrak6%
\font\fivefrak=cmfrak5%
\newfam\frakfam  \textfont\frakfam=\tenfrak  \scriptfont\frakfam=\sevenfrak  \scriptscriptfont\frakfam=\fivefrak%
\def\frak{\fam\frakfam\tenfrak}%

\font\tenmsam=msam10%
\font\eightmsam=msam8%
\font\sevenmsam=msam7%
\font\sixmsam=msam6%
\font\fivemsam=msam5%

\def\bb{\fam\bbfam\tenbb}%

\def\hexdigit#1{\ifnum#1<10 \number#1\else%
  \ifnum#1=10 A\else\ifnum#1=11 B\else\ifnum#1=12 C\else%
  \ifnum#1=13 D\else\ifnum#1=14 E\else\ifnum#1=15 F\fi%
  \fi\fi\fi\fi\fi\fi}
\newfam\msamfam  \textfont\msamfam=\tenmsam  \scriptfont\msamfam=\sevenmsam  \scriptscriptfont\msamfam=\fivemsam%
\def\msam{\msamfam\tenmsam}%
\mathchardef\leq"3\hexdigit\msamfam 36%
\mathchardef\geq"3\hexdigit\msamfam 3E%

\newdimen\bigsize  \bigsize=8.5pt
\newdimen\Bigsize  \Bigsize=11.5pt
\newdimen\biggsize \biggsize=14.5pt
\newdimen\Biggsize \Biggsize=17.5pt

\catcode`\@=11
\def\big#1{{\hbox{$\left#1\vbox to\bigsize{}\right.\n@space$}}}
\def\Big#1{{\hbox{$\left#1\vbox to\Bigsize{}\right.\n@space$}}}
\def\bigg#1{{\hbox{$\left#1\vbox to\biggsize{}\right.\n@space$}}}
\def\Bigg#1{{\hbox{$\left#1\vbox to\Biggsize{}\right.\n@space$}}}
\catcode`\@=12

\def\eightpoints{%

\def\rm{\fam0\eightrm}%
\textfont0=\eightrm   \scriptfont0=\sixrm   \scriptscriptfont0=\fiverm%
\textfont1=\eighti    \scriptfont1=\sixi    \scriptscriptfont1=\fivei%
\textfont2=\eightsy   \scriptfont2=\sixsy   \scriptscriptfont2=\fivesy%
\textfont3=\eightex   \scriptfont3=\sixex   \scriptscriptfont3=\fiveex%
\textfont\itfam=\eightit  \def\it{\fam\itfam\eightit}%
\textfont\slfam=\eightsl  \def\sl{\fam\slfam\eightsl}%
\textfont\ttfam=\eighttt  \def\tt{\fam\ttfam\eighttt}%
\textfont\bffam=\eightbf  \def\bf{\fam\bffam\eightbf}%

\textfont\frakfam=\eightfrak  \scriptfont\frakfam=\sixfrak \scriptscriptfont\frakfam=\fivefrak  \def\frak{\fam\frakfam\eightfrak}%
\textfont\bbfam=\eightbb      \scriptfont\bbfam=\sixbb     \scriptscriptfont\bbfam=\fivebb      \def\bb{\fam\bbfam\eightbb}%
\textfont\msamfam=\eightmsam  \scriptfont\msamfam=\sixmsam \scriptscriptfont\msamfam=\fivemsam  \def\msam{\msamfam\eightmsam}

\bigsize=6.8pt
\Bigsize=9.2pt
\biggsize=11.6pt
\Biggsize=14pt

\rm%
}

\def\anote#1#2#3{\smash{\kern#1in{\raise#2in\hbox{#3}}}%
  \nointerlineskip}     

\def\poorBold#1{\setbox1=\hbox{#1}\wd1=0pt\copy1\hskip.25pt\box1\hskip .25pt#1}

\mathchardef\lsim"3\hexdigit\msamfam 2E%
\mathchardef\gsim"3\hexdigit\msamfam 26%

\def\d{\,{\rm d}}
\def\D{\hbox{\rm D}}
\def\ds{\displaystyle}
\long\def\DoNotPrint#1{\relax}
\def\C{\hbox{\rm C}}
\def\E{\hbox{\rm E}}
 
\def\finetune#1{#1}
\def\fixedref#1{#1\note{fixedref$\{$#1$\}$}}
\def\g#1{g_{[0,#1)}}
\def\Gg{\Gamma(\gamma)}
\def\Ieps{{\bb I}_\epsilon}
\def\Id{\hbox{\rm Id}}
\def\K{{\rm K}}
\def\limB{\lim_{B\to\infty}}
\def\limeps{\lim_{\epsilon\to 0}}

\def\limsupn{\limsup_{n\to\infty}}
\def\limsupt{\limsup_{t\to\infty}}
\def\limn{\lim_{n\to\infty}}
\def\limt{\lim_{t\to\infty}}
\def\oF{\overline F{}}

\def\P{\hbox{\rm P}}
\def\qed{~~~{\vrule height .9ex width .8ex depth -.1ex}}
\def\sign{\hbox{\rm sign}}
\def\ss{\scriptstyle}
\def\St{{\cal S}_t}
\def\Var{\hbox{\rm Var}}

\def\calN{{\cal N}}
\def\calM{{\cal M}}

\def\calS{{\cal S}}

\def\II{\hbox{\tenbbm 1}}
\def\NN{{\bb N}\kern .5pt}
\def\RR{{\bb R}}
\def\ZZ{{\bb Z}}

\def\fS{{\frak S}}

\pageno=1

\centerline{\bf AN EXTENSION OF A LOGARITHMIC FORM}
\centerline{\bf OF CRAM\'ER'S RUIN THEOREM TO SOME FARIMA}
\centerline{\bf AND RELATED PROCESSES}

\bigskip

\centerline{Ph.\ Barbe$^{(1)}$ and W.P.\ McCormick$^{(2)}$}
\centerline{${}^{(1)}$CNRS, France, ${}^{(2)}$University of Georgia}

{\narrower
\baselineskip=9pt\parindent=0pt\eightpoints

\bigskip

{\bf Abstract.} 
Cram\'er's theorem provides an estimate for the tail probability of the 
maximum of a random walk with negative drift and increments having a moment
generating function finite in a neighborhood of the origin. The class of 
$(g,F)$-processes generalizes in a natural way random walks and fractional
ARIMA models used in time series analysis. For those $(g,F)$-processes with
negative drift, we obtain a logarithmic estimate of the tail probability
of their maximum, under conditions comparable to Cram\'er's. Furthermore,
we exhibit the most likely paths as well as the 
most likely behavior of the innovations leading to a large maximum.

\bigskip


\bigskip
 

}

\bigskip

\def\preveq{(\the\sectionnumber .\the\equanumber)}
\def\prevs{\the\sectionnumber .\the\snumber }

\section{Introduction.}
Cram\'er's theorem on the maximum of a random walk with negative drift provides
an estimate for the tail probability of this maximum when the moment 
generating function of the increments is finite in a neighborhood of the 
origin. Specifically, writing $M$ for the maximum of the random walk, it 
asserts that there are constants $c$ and $\theta$ such that
$$
  \P\{\, M>t\,\}\sim c e^{-\theta t}
  \eqno{\equa{cramer}}
$$
as $t$ tends to infinity; the constants $c$ and $\theta$ are explicit, but 
their formulas are irrelevant to the current discussion. We refer to 
Feller (1971, \S XI.7) for a proof of Cram\'er's theorem.

The purpose of this paper is to make a first step toward an extension of 
Cram\'er's result to a wider class of stochastic processes which encompass
some fractional ARIMA ones. As explained in Barbe and McCormick (2008) where
we dealt with the analogous problem in the heavy tail context, the motivations
are manifold. To summarize, besides the original application to insurance
mathematics which motivated Cram\'er, other areas of applications exist, 
such as queueing theory ---~the connection between risk and queueing theory
was pointed out in Prabhu (1961); see e.g.\ Janssen (1982) for an account
on this connection; furthermore, on a more fundamental level, a certain
analogy, described in Barbe and McCormick (2008),
has been developed between the asymptotic theory of the usual random
walk and that of some FARIMA processes, and it is natural to investigate to
which extent this analogy carries over in the context of Cram\'er's theorem.

Previous authors have considered ruin probabilities associated with
processes with dependent innovations. For instance, using a martingale
technique, Gerber (1982) considered bounded ARMA increments. His
result was extended by Promislow (1991) who removed the boundedness
assumption and dealt with a larger class of increments. In contrast,
using large deviations theory, building upon the work of Burton and 
Dehling (1990) as well as Iscoe, Ney and Nummelin (1985), 
Nyrhinen (1994, 1995, 1998) considered
increments following a stationary linear process with some having a
Markovian structure. M\"uller and Pflug (2001) extended some of
Nyrhinen's results by relating the asymptotic behavior of the moment
generating function of the ruin process at time $n$, as $n$ tends to
infinity, to the behavior of its maximum, hence, showing that the
G\"artner-Ellis (G\"artner, 1977; Ellis, 1984) approach in large
deviations leads to a ruin probability estimate. A common feature of
these works is that the processes under consideration exhibit short
range dependence in order to have an explicit behavior of some moment
generating functions. In contrast, the study of ruin probability
associated to continuous time processes has recently focussed on long
range dependent models. For instance, combining Duffield and O'Connell
(1995) with Chang, Yao and Zajic (1999) yields ruin probability
estimates for some nonnecessarily Gaussian long range memory processes
modeled after the fractional Brownian motion. More precise results
were obtained by H\"usler and Piterbarg (2004) for some Gaussian
processes. Our results may be viewed as a non-Gaussian and discrete
analogues of those continuous ones, in the sense that we are interested
in processes exhibiting long range dependence. Interestingly, for some
values of their paramaters, the processes considered in this paper,
suitably rescaled and normalized, converge to some fractional Brownian
motions.

A true extension of Cram\'er's theorem to FARIMA processes seems beyond what
one can achieve at the present, and we will only consider a logarithmic
form of it, namely, after taking the logarithm in \cramer,
$$
  \limt t^{-1}\log\P\{\, M>t\,\} = -\theta \, .
$$

The paper is organized as follows. The class of stochastic processes which we 
will consider and the main result are described in the next section. In
section 3  we describe the most likely scenario leading to a ruin, that is
to a large maximum of the processes under consideration. Section 4 contains
a broad outline of the proof. In section 5, we prove
some large deviations results which are of independent interest and lead
to the proof ---~inspired by Collamore (1996)~--- of the results of 
section 2. The result of section 3 is proved in section 6.

\bigskip

\noindent{\bf Notation.} Throughout this paper, if $(a_n)$ and $(b_n)$ are 
two sequences, we say that
'$a_n$ is lower bounded from above by an equivalent of $b_n$' and write
$a_n\lsim b_n$ if $a_n\leq b_n\bigl(1+o(1)\bigr)$ as $n$ tends to infinity,
or, equivalently, if $\limsup_{n\to\infty} a_n/b_n\leq 1$. The symbol $\gsim$
is defined in an analogous way.

\bigskip


\section{Main result.}
Barbe and McCormick (2008) introduced $(g,F)$-processes as a natural extension
of FARIMA processes. To define such a process, we start with a 
function $g$ which
is real analytic on $(-1,1)$ and a distribution function $F$. The function $g$
has a Taylor series expansion
$$
  g(x)=\sum_{i\geq 0} g_i x^i \, .
$$
Considering a sequence $(X_i)_{i\geq 1}$ of independent and identically 
distributed random variables with common distribution function $F$, we define
the $(g,F)$-process $(S_n)_{n\geq 0}$ by $S_0=0$ and
$$
  S_n=\sum_{0\leq i<n} g_i X_{n-i} \, .
$$
When $g(x)=(1-x)^{-1}$, the corresponding process is the random walk associated
to the sequence $(X_n)_{n\geq 1}$. Some nonstationary ARMA processes are 
obtained when $g$ is a rational function, and FARIMA processes are obtained
when $g(x)$ is the product of some negative power of $1-x$ and a rational 
function in $x$.

For the process to drift toward minus infinity and mimic the behavior of the
random walk involved in Cram\'er's theorem, it is natural to impose that the
mean $\mu$ of $F$ is negative and that
$$
  \limn \sum_{0\leq i<n} g_i = +\infty \, .
  \eqno{\equa{gdrift}}
$$
Indeed, in this case, the expectation of $S_n$ drifts toward minus infinity. 
A consequence of \gdrift\ is that $g$ has a singularity at $1$. To obtain a
satisfactory theory, we need to restrict the type of singularity by assuming
that $g$ is regularly varying at $1$ of positive index, meaning, as explained 
for instance in
Bingham, Goldie and Teugels (1989), that there exists a positive $\gamma$
such that for any positive $\lambda$,
$$
  \limt {g(1-1/\lambda t)\over g(1-1/t)} = \lambda^\gamma \, .
$$
This assumption is satisfied by ARIMA processes.

Let $\Id$ be the identity function on the real line.
We then consider a function $U$, defined up to asymptotic equivalence by
the requirement
$$
  g(1-1/U)\sim \Id
$$
at infinity. This function, which plays a key role in our result, is 
necessarily regularly varying at infinity of index $1/\gamma$. However,
for notational simplicity, writing $\Gamma(\cdot)$ for the gamma function, 
it will be better to use the function
$$
  V=\Gamma(1+\gamma)^{1/\gamma} U \, ,
$$
which could alternatively be defined by the requirement $g(1-1/V)\sim 
\Gamma(1+\gamma)\Id$ at infinity.

In order to concentrate on the principles and the key arguments, we assume
throughout this paper that the coefficients $g_i$ are nonnegative. This
restriction can be overcome with the introduction of the proper tail balance
condition.

To have a compact notation, we introduce the kernel
$$
  k_\gamma(u)=\cases{ \gamma (1-u)^{\gamma-1} & if $0\leq u<1$, \cr
                      \noalign{\vskip 3pt}
                      0 & if $u\geq 1$,\cr }
$$
defined on the nonnegative half-line.

Further notation related to large deviations theory is needed in order to 
state our main result. As the proof shows, the appearance of some large 
deviations formalism is not coincidental. Cram\'er's theorem assumes
that the moment generating function
$$
  \varphi(\lambda)=\E e^{\lambda X_1}
$$
is finite in a neighborhood of the origin. A classical consequence of
H\"older's inequality is that $\log\varphi$ is convex. This implies that the
function
$$
  \lambda\mapsto \int_0^1\log\varphi \bigl(\lambda k_\gamma(u)\bigr)\d u
  \eqno{\equa{intMGF}}
$$
is convex as well on its domain. This function will be of importance in 
our results. It is not clear a priori that this function is nontrivial
in the sense that if $\gamma$ is less than $1$ it could be infinite for all
nonvanishing $\lambda$. This suggests that we should consider two cases, 
according to the finiteness of the integral involved in \intMGF.

The convex conjugate (see e.g.\ Rockafellar, 1970) of 
the function involved in \intMGF, at a nonnegative argument $a$, is
$$
  J(a)=\sup_{\lambda>0}\,
  \Bigl( a\lambda 
         -\int_0^1\log\varphi\bigl(\lambda k_\gamma(u)\bigr)\d u
  \Bigr)
  \, .
$$

To a moment generating function $\varphi$ one also associates the 
corresponding mean function $m$, which is the derivative
$(\log\varphi)'$ ---~see Barndorff-Nielsen (1978), Brown (1986) or
Letac (1992).

The following convention will be convenient. We say that a $(g,F)$-process
satisfies the standard assumption if it satisfies the following

\bigskip

\noindent{\bf Standard assumption.}\quad {\it 
  The function $g$ is regularly varying
  of positive index at $1$ and its coefficients $(g_i)_{i\geq 0}$ are 
  nonnegative. Moreover $g_0$ does not vanish. In case the sequence
  $(g_n)_{n\geq 0}$ converges to $0$, it is asymptotically equivalent
  to a monotone sequence.
  The distribution function $F$ has a moment generating function finite
  on the nonnegative half-line. The image of the mean function contains
  the half line $[\, 0,\infty)$.}

\bigskip

With respect to the monotonicity requirement for the sequence 
$(g_n)_{n\geq 0}$ involved in the standard assumption,
it will follow from Proposition 1.5.3 in Bingham, Goldie and Teugels (1989) 
and Lemma \fixedref{5.1.1} that regular variation of $g$ implies that
$(g_n)_{n\geq 1}$ is asymptotically equivalent to a monotone sequence
whenever the index of regular variation of $g$ is different from $1$.

Let $(S_n)_{n\geq 0}$ be a $(g,F)$-process. If the first $k$ coefficients
$g_0$, $g_1, \ldots , g_{k-1}$ vanish and $g_k$ does not, then 
$(S_{n+k})_{n\geq 1}$ is a $(g/\Id^k,F)$-process, and the first Taylor 
coeffecient of $g/\Id^k$ does not vanish. Thus, in the standard assumption,
the condition that $g_0$ does not vanish bears no restriction.

Note that in the standard assumption, the condition on the moment generating
function is stronger than in Cram\'er's theorem. 
The assumption on the mean function is a rather standard one in large 
deviations theory. H\"older's inequality implies that $\log\varphi$ is convex
and the mean function is nondecreasing. Our assumption ensures that the
equation $m(\lambda)=x$ has a solution for every positive $x$.

We also say that a $(g,F)$-process satisfying the standard assumption
has a negative mean if its expectation is negative at all time. Since the
innovations are independent and identically distributed, considering the 
expectation of the process at time $1$, this is equivalent to require that
the innovations have negative mean.

Our first result treats the case where the integral \intMGF\ is finite.
It calls for many remarks, stated after the theorem, which will clarify 
both the hypotheses and the conclusion.

\Theorem{\label{ThCaseA}
  Consider a negative mean $(g,F)$-process which satisfies the standard 
  assumption. Assume that either one of the following conditions hold:
  
  \smallskip

  \noindent (i) $\limsup_{n\to\infty}\max_{0\leq i<n} g_i/g_n$ is finite;

  \smallskip

  \noindent (ii) $\limn \max_{0\leq i\leq n} g_i/g_n=+\infty$ 
  and $-\log\oF$ is regularly varying with
  index $\alpha$ such that $\alpha\gamma>1$; moreover, $m'$ is regularly
  varying.

  \smallskip

  \noindent Then, the function $J$ is defined and finite on the nonnegative 
  half-line and the maximum $M$ of the $(g,F)$-process satisfies
  $$
    \limt  V(t)^{-1} \log \P\{\, M>t\,\}
    = -\inf_{x>0} x J(x^{-\gamma}) \, .
  $$
}

We now make some remarks on the conclusion of the theorem, which will be
followed by remarks on its assumptions.

Writing $\theta$ for the negative of the limit involved in its statement, 
this theorem asserts that
$$
  \P\{\, M>t\,\}\sim e^{-\theta V(t)(1+o(1))}
$$
as $t$ tends to infinity. This leads to the following observation which may 
constitute a caveat of pedagogical value. Fix the distribution function $F$
and consider the analytic function $g$ as a parameter. As we increase its
singularity at $1$, the process drifts toward minus infinity at a faster rate,
for its mean at time $n$ is $\mu\sum_{0\leq i<n}g_i$. One might guess that this
makes it harder for the process to reach a high threshold. However, our theorem
asserts that the logarithmic order of this probability is $-V(t)$, which 
becomes larger with $g$. So, making the mean to diverge to minus infinity 
faster, makes it more likely for the process to reach a high level! 
This phenomenon will be explained in the next section.

In the same spirit, it will follow from equality \fixedref{(5.2.13)}
that multiplying the $X_i$ by a scale factor $\sigma$ divides $\theta$ by
$\sigma^{1/\gamma}$. Thus, increasing the drift toward minus infinity through 
a scaling increases the likelihood for $M$ to 
take very large values.

On a different note, we see that as in Cram\'er's theorem, the tail of
the distribution function of the increments is involved in the
conclusion of Theorem \ThCaseA\ only in the constant $\theta$ and not
in the logarithmic decay $V$.

It is also of interest to note that if $\gamma$ is greater than $1$, then
$V\ll\Id$ at infinity. In this case, Theorem \ThCaseA\ shows that the
distribution of the maximum of the process is subexponential, even
though the innovations are superexponential. Such a possibility was
observed in a different context by Kesten (1973).

Regarding the assumptions of Theorem \ThCaseA, note that in case (i)
we must have $\gamma$ at least $1$. In case (ii), the condition
that $\max_{0\leq i\leq n} g_i/g_n$ diverges to infinity is 
equivalent to the convergence of $(g_n)_{n\geq 0}$ to $0$, which
forces $\gamma$ to be at most $1$.

Let $\beta$ be the conjugate exponent of $\alpha$, that is such that
$\alpha^{-1}+\beta^{-1}=1$. It follows from Kasahara's theorem
(Bingham, Goldie and Teugels, 1989, Theorem 4.12.7) that $-\log\oF$ is
regularly varying of index $\alpha$ if and only if $\log\varphi$ is
regularly varying of index $\beta$. Since $\log\varphi$ is convex, its
derivative is monotone, and the monotone density theorem combined with
Kasahara's theorem implies that $-\log\oF$ is regularly varying of
index $\alpha$ if and only if $m$ is regularly varying of index
$\beta-1$. The assumption on $m'$ in Theorem \ThCaseA\ is stronger.
This assumption is not completely satisfactory since its meaning in
terms of the distribution function is not clear.

Under the assumptions of Theorem \ThCaseA, we must have
$\beta-2>-1$. Hence, using Karamata's theorem in addition to the
previous paragraph, we see that the assumption of Theorem \ThCaseA\ on
$-\log \oF$ and $m'$ is equivalent to the single assumption that $m'$
is regularly varying of index $\beta-2$ with $\beta(1-\gamma)<1$.

\bigskip

Our second result considers the case where the integral involved in \intMGF\
is infinite, and hence the function $J$ in Theorem \ThCaseA\ is not
defined. This essentially occurs when 
$\gamma$ is less than $1$ and $\alpha\gamma$ is at most $1$. If $\gamma$ is
less than $1/2$ then the centered process $S_n-\E S_n$ converges in 
$\hbox{\rm L}^2$. For $\gamma$ less than $1/2$, let $Z_n$ be the 
linear process
$$
  Z_n=\sum_{i\geq 0} g_i (X_{n-i}-\mu) \, .
$$
In this case, we see that the ruin  problem for $S_n$ is rather similar to that
of determining the probability that the process $(Z_n)_{n\geq 1}$ crosses
the moving boundary $t-ES_n$. This problem is of somewhat different nature
than what is the focus of this paper, for the centered process is well
approximated by a stationary one. Therefore, we will limit ourselves
to the case where $\gamma$ is greater than $1/2$.

We write $|g|_\beta$ for
the $\ell_\beta$-norm of the sequence of its coefficients, that is for
$\bigl(\sum_{i\geq 0} g_i^\beta\bigr)^{1/\beta}$.

\Theorem{\label{ThCaseC}
  Consider a $(g,F)$-process which satisfies the standard assumption
  and with $1/2<\gamma<1$.
  Assume furthermore that $-\log\oF$ is regularly 
  varying of index $\alpha$ greater than $1$ and that $\alpha\gamma<1$. 
  Let $\beta$ be the conjugate exponent
  of $\alpha$. Then, the maximum $M$ of the process satisfies
  $$
    \log \P\{\, M>t\,\}
    \sim |g|_\beta^{-\alpha} \log\oF(t) \, .
  $$
  as $t$ tends to infinity.
}

\bigskip

Comparing Theorems \ThCaseC\ and \ThCaseA, we see that
in Theorem \ThCaseC, the condition $\alpha\gamma<1$ forces the rate of 
growth of $-\log F(t)$, regularly varying of index $\alpha$, to be much 
slower than that of $U(t)$, regularly varying of index $1/\gamma$. 

\bigskip


\section{How to go bankrupt?}
The purpose of this section is to determine the most likely paths which
lead to the maximum of our $(g,F)$-processes to reach a high
threshold. Beyond its relevance to choosing interesting alternatives in
change point problems, in the context of ruin probability, this amounts to find
the most likely way of becoming bankrupted. In a
different context, high risk scenarios have been the subject of
Balkema and Embrechts (2007) monograph where further discussion of the
topic may be found. More closely related to the topic of this paper,
is the work of Chang, Yao and Zajic (1999) in the continuous setting,
who consider the analogous problem for fractional integrals of
continuous time processes. In fact we are seeking more information.
Not only are we interested in the most likely paths, but we would also like
to understand how they arise, and, therefore, have a description of
the innovations as well. In the heavy tail case, it is shown in Barbe
and McCormick (2008) that a large value of the maximum of the process
is most likely caused by a large value of an innovation. In contrast,
in a slightly different setting than that of the current paper, but
nonetheless related, for the usual random walk, Csisz\'ar (1984,
Theorem 1) shows that a large deviation is likely caused by a
cooperative behavior of the increments which pushes the sum
upward. More precisely, Csisz\'ar's result implies that the
conditional distribution of the first increment, given that the sum
$S_n$ exceeds an unlikely threshold $nu$, converges to the
distribution $\d F_u(x)=e^{m^\leftarrow (u)x} \d F(x)/\varphi\circ
m^\leftarrow (u)$. The distribution $\d F_u$ has mean $u$. For the
usual random walk, since the increments are exchangeable given their
sum, Csisz\'ar's result asserts that, loosely, a randomly chosen
increment, or a typical increment, has a conditional distribution
about $\d F_u$.  Thus, asymptotically, the bulk of the increments
behave like a random variable of mean $u$ under the conditional
distribution that the random walk at time $n$ exceeds $nu$. We refer
to Diaconis and Freedman (1988) for a refined result in the framework of
exponential families.

In general, for $(g,F)$-processes, the innovations are not exchangeable
given the value of the process at time $n$, and, paralleling what has 
been done for the random walk, it is of interest to identify the 
cooperative behavior of the increments, if any, which makes the process 
to reach a high level.

Besides a theoretical understanding, this type of conditional limiting
result has some bearing on simulation techniques of rare events by importance
sampling (Hammersley and Handscomb, 1964). Indeed, when specialized to the 
regular random walk, Sadowsky (1996) gives a rationale for using the
limiting conditional distribution of the increments to simulate unlikely
paths of random walks using importance sampling; see also Dieker and Mandjes 
(2006). Our result is a key building
block to extend this technique to some FARIMA processes, and, more generally,
to $(g,F)$-processes.

To investigate these questions, we consider first the rescaled trajectory
$$
  \St (\lambda)=S_{\lfloor \lambda V(t)\rfloor}/t \, , \qquad 
  \lambda \geq 0 \, .
$$
Next, to study the behavior of the innovation, we consider the sequential
measure
$$
  \calM_t={1\over V(t)}\sum_{i\geq 1} \delta_{(i/V(t),X_i)}
$$
which puts mass $1/V(t)$ at each pair $\bigl(i/V(t),X_i\bigr)$. In
contrast with a standard empirical measure which would put equal mass
on each innovation up to some fixed time, the sequential measure keeps
track of the sequential ordering of the innovation through the first
component $i/V(t)$.

Of further interest is also the normalized first time that the process 
reaches the level $t$,
$$
  \calN_t={1\over V(t)}\min\{\, n\, :\, S_n>t\,\} \, .
$$
In order to speak of convergence of the stochastic process $\St$, we view it 
in the Skorohod space $\D[\,0,\infty)$ equipped with the Skorohod topology
(Billingsley, 1968; Lindvall, 1973). 

In what follows, we call $[\,0,\infty)\times \RR$ the right half-space. 
A subset of the right half-space of the form $[\,a,b\,]\times \RR$ is called a
vertical strip.

The measure $\calM_t$ belongs to the space $\calM([\,0,\infty)\times
\RR)$ of $\sigma$-finite measures on the right half-space. We consider
this space equipped with a topology between those of vague and weak$*$
convergences defined as follows.  Let 
$\C_{\K,{\rm b}}([\,0,\infty)\times \RR)$ be the space of all real-valued
continuous and bounded functions on the right half-space, supported on
a vertical strip. A basis for the topology on
$\calM([\,0,\infty)\times \RR)$ is defined by the sets
$$
  \Bigl\{\, \mu\in\calM([\,0,\infty)\times \RR) \, :\, \forall i=1,\ldots ,k
  \, ,\,
  \Bigl|\int f_i \d (\mu-\nu)\Bigr|<\epsilon\,\Bigr\} \, ,
$$
indexed by
$$
  \nu\in\calM([\,0,\infty)\times \RR) \, ,\quad 
  f_i\in \C_{\K,{\rm b}}([\,0,\infty)\times \RR) \, ,\quad 
  \epsilon>0 \, .
$$
In this paper, except specified otherwise, all convergences of measures
on the right half-space are for this topology.

Our next result gives the limit in probability of the various quantities
introduced, conditionally on having the process reaching the level $t$, and 
under the assumptions of Theorem \ThCaseA. We assume that
$$
  \tau=\arg\min_{x>0} xJ(x^{-\gamma})
  \hbox{ is unique.}
  \eqno{\equa{HypTauUnique}}
$$
Furthermore, we define the constant $A$ to be the solution of
$$
  \tau^{-\gamma}
  =\int_0^1 k_\gamma (u) m\bigl(Ak_\gamma(u)\bigr)\d u \, .
  \eqno{\equa{ADefA}}
$$
Let $L$ be the Lebesgue measure. We define the measure $\calM$ by its 
density with respect to the product measure $L\otimes F$,
$$
  {\d\calM\over \d (L\otimes F)}(v,x)
  = { \exp \bigl( Ak_\gamma(v/\tau)x\bigr)
      \over 
      \varphi\bigl( Ak_\gamma(v/\tau)\bigr) }
  \, .
  \eqno{\equa{MDef}}
$$
In particular, since $k_\gamma$ vanishes on $[\,1,\infty)$, the measure
$\calM$ coincides with $L\otimes F$ on $[\,\tau,\infty)\times\RR$.
We also define the function
$$
  \calS (\lambda)
  =\int_0^\lambda \gamma (\lambda-v)^{\gamma-1}
    m\bigl(Ak_\gamma(v/\tau)\bigr) \d v \, . 
  \eqno{\equa{SDef}}
$$
Writing
$$
  \calS(\lambda)
  =\lambda^\gamma\int_0^1 k_\gamma (v) 
    m\bigl(Ak_\gamma(v\lambda/\tau)\bigr) \d v
$$
and using \ADefA, we see that $\calS(\tau)=1$. 

The following result describes the most likely ruin scenario.

\Theorem{\label{ThConditionalA}
  Under the assumptions of Theorem \ThCaseA, the following hold in probability
  under the conditional probability given $M>t$ as $t$ tends to infinity:

  \smallskip

  \noindent (i) $\calN_t$ converges to $\tau$;

  \smallskip

  \noindent (ii) $\calM_t$ converges to $\calM$;

  \smallskip

  \noindent (iii) moreover, if the moment generating function of $|X_1|$ is
  finite in a neighborhood of the origin, then $\St$ converges locally 
  uniformly to $\calS$.
}

\bigskip

Regarding the hypotheses of Theorem \ThConditionalA, under those of
Theorem \ThCaseA,
the assumption that the moment generating function of $|X_1|$ is finite
in a neighborhood of the origin
is weaker than a tail balance condition. A close look at the proof shows
that this assumption is used only to prove assertion (ii) ---~see Lemma 
\fixedref{6.3.1}.

\bigskip

Loosely speaking, the meaning of assertion (ii) is that the conditional
distribution of $X_{\lfloor vV(t)\rfloor}$ given $M>t$ converges to the
measure
$$
  \cases{ { \ds e^{A\gamma(1-v/\tau)^{\gamma-1} x}
            \over\ds \varphi \bigl(A\gamma(1-v/\tau)^{\gamma-1}\bigr) } \d F(x)
          & if $v\leq \tau$ \cr\noalign{\vskip 3pt}
            \d F(x) 
          & if $v> \tau$,\cr 
        }
$$
with mean $m\bigl( A\gamma(1-v/\tau)^{\gamma-1}\bigr)$ if $v<\tau$, and
$\mu$ if $v\geq\tau$.
Thus it asserts that a large value of $M$ is likely caused by a cooperative
behavior of the random variables up to a time $\tau V(t) \bigl(1+o(1)\bigr)$,
while the remainder of the innovations keep their original 
distribution. This somewhat confirms that the Cram\'er ruin model might be 
unrealistic in some situations. Indeed, Theorem \ThConditionalA\ shows that
for the process to reach the large level $t$, both the increments and the
process, from the very beginning, have to follow a very unlikely
path. One would think that seeing such a strange path unfolding, a careful
insurer would quickly reexamine the model and raise the premium accordingly.

Theorem \ThConditionalA\ also explains why Theorem \ThCaseA\ implies that
for those $(g,F)$-processes, adding more drift toward minus infinity may
increase the likelihood of a large maximum. Indeed Theorem \ThConditionalA\
indicates that a large value of the maximum is likely to be caused by
many innovation being large; but if the weights $(g_n)_{n\geq 0}$ are
made larger, then comparatively smaller innovation suffices for the maximum
of the process to reach a large value, because the coefficients 
$(g_n)_{n\geq 0}$ amplify the innovations.

\medskip

We now consider an example of processes of interest and for which the 
limit involved in Theorems \ThCaseA\ or \ThCaseC\ can be made explicit.
In general this limit must be evaluated by numerical methods.

We consider a Gaussian FARIMA process. More specifically, we consider $F$
to be the Gaussian distribution function with mean $\mu$ and variance
$\sigma^2$, and we introduce two 
polynomials $\Theta$ and $\Phi$, neither of which vanishes at $1$.
We consider
the function $g(x)=(1-x)^{-\gamma}\Theta(x)/\Phi(x)$, so that the corresponding
$(g,F)$-process is a FARIMA($\Phi,\gamma,\Theta$) process whose innovations 
have a common distribution function $F$. For this specific function $g$ we may
take $U(t)=\bigl(t\Phi(1)/\Theta(1)\bigr)^{1/\gamma}$.

The moment generating function of the innovations is
$$
  \varphi(\lambda)=e^{\lambda\mu+\sigma^2\lambda^2/2} \, .
$$
The function involved in \intMGF\ is then
$$
  \gamma\int_0^1 \mu\lambda\gamma u^{\gamma-1} + {\sigma^2\over 2}
  (\lambda\gamma u^{\gamma-1})^2 \d u
  = \lambda\mu + {\sigma^2\over 2}\lambda^2 {\gamma^2\over 2\gamma-1} \, .
$$
This implies that
$$\eqalign{
  J(a)
  &{}=\sup_\lambda
   \Bigl( a\lambda -\lambda\mu-{\sigma^2\over 2} 
   \lambda^2 {\gamma^2\over 2\gamma-1}\Bigr)
  \cr
  &{}={(a-\mu)^2 (2\gamma-1)\over 2\sigma^2\gamma^2} \, . \cr
  }
$$
Using standard calculus one more time, we obtain
$$
  \inf_{x>0}xJ(x^{-\gamma}) 
  = {2(2\gamma-1)^{1/\gamma-1}\over\sigma^2}(-\mu)^{2-1/\gamma} \, .
$$
Therefore, the conclusion of Theorem \ThCaseA\ is that
$$\displaylines{\qquad
  \log\P\{\, M> t\,\}
  \hfill\cr\hfill
  \sim -t^{1/\gamma}\Bigl({\Phi(1)\over\Theta(1)}\Bigr)^{1/\gamma}
  \Gamma(1+\gamma)^{1/\gamma}
  2(2\gamma-1)^{(1/\gamma)-1} \Bigl({\mu\over\sigma}\Bigr)^2 (-\mu)^{-1/\gamma}
  \qquad\cr}
$$
as $t$ tends to infinity.

To calculate the limiting process $\calS$, for simplicity we restrict ourselves
to the case where the mean $\mu$ is $-1$ and the standard deviation $\sigma$
is $1$. Then, $m(\lambda)=\lambda-1$, and
$$
  \calS(\lambda)
  = A\int_0^{\lambda\wedge\tau} \gamma (\lambda-v)^{\gamma-1} 
    \gamma\Bigl(1-{v\over\tau}\Bigr)^{\gamma-1} \d v -\lambda^\gamma \, .
$$

The following graphic shows the shape of the limiting function when $\gamma$
is $2/3$, $1$ and $2$.

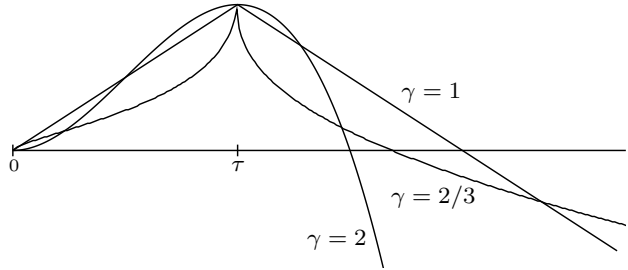
\setbox1=\hbox to 3.2truein{\eightpoints\vbox to 1.4in{
  \kern 2.1in
  \includegraphics{logCramer.fig}
  \anote{-.02}{1.15}{$\ss 0$}
  \anote{1.04}{1.15}{$\tau$}
  \anote{1.4}{.8}{$\gamma=2$}
  \anote{1.8}{1}{$\gamma=2/3$}
  \anote{1.85}{1.5}{$\gamma=1$}
  \vfill
  \vss
}\hfill}    

\bigskip

\centerline{\box1}

\bigskip

We conclude this section by some remarks concerning Theorem
\ThConditionalA\ and its proof. A close look at the proofs of Theorems
\ThCaseA\ and \ThConditionalA\ reveals that the same technique allows one
to derive a large deviations principle for the process $\calS_t$ and
the measure $\calM_t$ under the conditional distribution of $M$
exceeding $t$, as $t$ tends to infinity, in the spirit of Collamore
(1998).

One can also see that under the assumptions of Theorem
\ThCaseA, the large deviations principle for FARIMA processes proved in
Barbe and Broniatowski (1998) remains true when the order of
differentiation $\gamma$ is between $1/2$ and $1$ and that the
logarithm of the tail of the distribution function of the innovation is
regularly varying with index greater than $1/\gamma$. This has the
following interesting consequence about the standard partial sum
process, $\Pi_n(\lambda)=n^{-1}\sum_{1\leq i\leq n\lambda} X_i$,
$0\leq\lambda\leq 1$. Consider the Cram\'er transform of the
increment, $I(x)=\sup_\lambda \lambda x-\log\varphi(\lambda)$.
Mogulskii's (1976) theorem (see also Dembo and Zeitouni, 1993, \S 5.1)
asserts that the partial sum process obeys a large deviations
principle, in the supremum norm topology. We can write 
the partial sum process at time $t$
as $\int_0^\lambda \d \Pi_n(v) =\int \II_{[0,\lambda)}(v)\d \Pi_n(v)$.
One could then wonder if some fractional integral of $\Pi_n$ still
obeys a large deviations principle. While an integration by parts shows
that for $\gamma$ greater than $1$, the process $\lambda\in
[\,0,1\,]\mapsto \int_0^\lambda (\lambda-v)^{\gamma-1}\d\Pi_n(v)$ obeys a
large deviation, the proof of Theorem \ThConditionalA\ shows that such
a large deviations principle still holds if $1/2<\gamma<1$, provided
that $\log\oF$ is regularly varying of index greater than
$1/\gamma$. The Gaussian case, $\alpha=2$ appears to be a boundary one
corresponding to $\gamma=1/2$; and this matches the fact that the
Brownian motion belongs to any set of functions with H\"older exponent
less than $1/2$.

\bigskip


\section{Generalities.}
The study of first passage times using large deviations is by now a classical
topic which has been presented in book form by
Freidlin and Wentzell's (1984).
The purpose of this section is to give another short variation
on this theme, with a formalism more suitable for the problems 
considered in this paper. What follows is inspired by the work of Collamore
(1998) as well as Duffield and Whitt (1998). However, in contrast to those
authors, we are interested in processes which are not Markovian, not mixing
and not monotone.

Some notation will purposely be identical to those used in the previous
sections, the reason being that they have the same meaning when specialized
to the context of the previous sections; this will be clear during the proofs
of Theorems \ThCaseA, \ThCaseC\ and \ThConditionalA.

In what follows, sequences are viewed as functions defined on 
the nonnegative half-line and evaluated at the integers. Therefore, if we write
$(a_n)_{n\geq 1}$ for a sequence, we will also speak of the function $a$,
meaning that $a_n=a(n)$ for every positive integer $n$. If we are given the
sequence, it is understood that the function $a$ is obtained by a linear
interpolation say; other `reasonable' interpolation procedures would do just
as well.

In this section we consider a stochastic process $(S^0_n)_{n\geq 1}$ and a
sequence $(s_n)_{n\geq 1}$ which diverges to infinity. We are 
interested in evaluating the probability that the process $(S^0_n)_{n\geq 1}$
crosses the moving boundary $(t+s_n)_{n\geq 1}$ for large values $t$.
In other words, assuming that $M=\max_{n\geq 1}S^0_n-s_n$ is well defined,
we are interested in finding an estimate of
$$
  \P\{\, \exists n\geq 1\,:\, S^0_n>t+s_n\,\}
  = \P\{\, M>t\,\}
$$
as $t$ tends to infinity. Assuming that the function
$$
  \hbox{$s$ is regularly varying of positive index $\gamma$,}
  \eqno{\equa{HypSRV}}
$$
there exists a function $V$, defined, up to asymptotic equivalence, by
the relation $s\circ V\sim \Id$ at infinity.  Also of interest is the
normalized first passage time at which the process crosses the moving
boundary,
$$
  \calN_t={1\over V(t)}\min\{\, n\geq 1\, : \, S^0_n>t+s_n\,\} \, .
$$

Suppose that $(S^0_n)_{n\geq 1}$ obeys a large deviations principle in the
sense that there exist two functions $r$ and $I$ such that for any positive
$x$
$$
  \log\P\{\, S^0_n> s_nx\,\}\sim -r_n I(x)
  \eqno{\equa{HypLDP}}
$$
as $n$ tends to infinity. Since the left hand side of \HypLDP\ is montone
in $x$, so is the right hand side, and, necessarily, $I$ is monotone as well
as continuous almost everywhere. If we assume more, namely that
$$
  \hbox{$I$ is continuous,}
  \eqno{\equa{HypICont}}
$$
then the asymptotic equivalence in \HypLDP\ holds locally uniformly in
$x$ over the nonnegative half-line, because a pointwise convergent
sequence of nondecreasing functions whose limit is continuous
converges locally uniformly (see Rudin, 1976, chapter 7, exercise 13).

For our problem, we will be able to assume that
$$
  \hbox{$r$ is regularly varying of positive index $\rho$.}
  \eqno{\equa{HypRRv}}
$$
In this case, $r$ is asymptotically equivalent to a nondecreasing function,
and we will consider, without any loss of generality, that $r$ is 
nondecreasing. We define $\theta$ as
$$
  \theta=\inf_{x>0} x^\rho I(x^{-\gamma}+1)\, .
  \eqno{\equa{ThetaDef}}
$$

We will also assume that the process is unlikely to reach the moving 
boundary $t+s_n$ before a time of order $V(t)$, in the sense that
$$
  \limeps\limsupt {1\over r\circ V(t)} 
  \log\P\{\, \exists n \,:\, 1\leq n\leq\epsilon V(t) \, ; \, S^0_n>t+s_n\,\} 
  \leq -\theta \, .
  \eqno{\equa{HypSmallN}}
$$

Equipped with these perhaps drastically looking ---~but to be proved useful~---
conditions, we have the following.

\Proposition{\label{MasterTh}
  If \HypSRV--\HypSmallN\ hold, then
  $$
    \log\P\{\, M>t\,\}\sim -\theta r\circ V(t)
  $$
  as $t$ tends to infinity. Moroever, if
  $$
  \tau=\arg\min_{x>0} x^\rho I(x^{-\gamma}+1)\hbox{ exists and is unique,}
  \eqno{\equa{HypMinUnique}}
  $$
  then $\calN_t$ converges to $\tau$ in probability given $M>t$, as $t$
  tends to infinity.
  }

\bigskip

\noindent{\bf Remark.} If we replace assumption \HypLDP\ by
$$
  \log\P\{\, S^0_n>x s_n\,\}\lsim -r_n I(x)
  \eqno{\equa{HypUBLDP}}
$$
as $n$ tends to infinity, the proof of Proposition \MasterTh\ shows that
$$
  \log\P\{\, M>t\,\} \lsim -\theta r\circ V(t)
$$
as $t$ tends to infinity. This remark will be useful to prove Theorem \ThCaseC.

\bigskip

In order to prove Proposition \MasterTh, we need the following lemma.

\Lemma{\label{LemmaMasterTh}
  Let $r$ be a nondecreasing regularly varying function of positive index.
  Then
  $$
    \sum_{n\geq k}e^{-r_k}\lsim k e^{-r_k}
  $$
  as $k$ tends to infinity.
}

\bigskip

\Proof Since $r$ is nondecreasing,
$$\eqalignno{
  \sum_{n\geq k+1} e^{-r_n}
  &{}\leq \int_k^\infty e^{-r(x)} \d x \cr
  &{}= e^{-r_k} \int \II\{\, k\leq x\, ;\, r(x)\leq r(k)+u\,\} e^{-u}\d u\,
    \d x\, .\qquad
  &\equa{EqLemmaMasterThA} \cr
  }
$$
Let $\rho$ be the index of regular variation of $r$. Let $\epsilon$ be an
arbitrary positive real number. Using Potter's bound, if $k$ is large enough
and $x\geq k$ then $r(x)/r(k)\geq (1-\epsilon)(x/k)^{\rho-\epsilon}$. In 
particular, if moreover $r(x)\leq r(k)+u$, then
$$
  x\leq k\Bigl( {1\over 1-\epsilon}\Bigl( 1+{u\over r(k)}\Bigr)
         \Bigr)^{1/(\rho-\epsilon)} \, .
$$
Thus, for $k$ large enough and after permuting the integration with respect
to $u$ and $x$, \EqLemmaMasterThA\ is at most
$$
  k{e^{-r_k}\over (1-\epsilon)^{1/(\rho-\epsilon)}}
  \int_0^\infty \Bigl( 1+{u\over r(k)}\Bigr)^{1/(\rho-\epsilon)} e^{-u} \d u
  \, .
$$
It follows from Lebesgue's dominated convergence theorem that this bound is
asymptotically equivalent to $ke^{-r_k}/(1-\epsilon)^{1/(\rho-\epsilon)}$
as $k$ tends to infinity. Since $\epsilon$ is arbitrary, this yields the 
result.\hfill\qed

\bigskip

\noindent{\bf Proof of Proposition \MasterTh.} The proof of the first
assertion consists in establishing the proper upper and lower bounds.

\noindent{\it Upper bound.} Let $\epsilon$ be a positive real number less than
$1$. For $t$ large enough and uniformly in $n$ between $\epsilon V(t)$ and
$V(t)/\epsilon$,
$$
  r_n
  =r\Bigl( V(t){n\over V(t)}\Bigr) 
  \sim r\circ V(t) \Bigl( {n\over V(t)}\Bigr)^\rho
$$
and
$$
  {t\over s_n}
  = {t\over s\Bigl( V(t){\ds n\over\ds V(t)}\Bigr)}
  \sim \Bigl( {n\over V(t)}\Bigr)^{-\gamma}
$$
as $t$ tends to infinity. In particular,
$$\eqalign{
  r_nI\Bigl( {t\over s_n}+1\Bigr) 
  &{}\sim r\circ V(t) \Bigl({n\over V(t)}\Bigr)^\rho 
    I\Bigl( \Bigl( {n\over V(t)}\Bigr)^{-\gamma}+1\Bigr) \cr
  &{}\gsim r\circ V(t)\theta \, . \cr
  }
$$
Combining this lower bound with the large deviations assumption
\HypLDP\ yields, in the range of $n$ between $\epsilon V(t)$ and
$V(t)/\epsilon$ and for $t$ large enough,
$$
  \P\{\, S^0_n>t+s_n\,\}
  \leq \exp\bigl( -r\circ V(t) \theta (1-\epsilon)\bigr) \, . 
  \eqno{\equa{EqMasterThA}}
$$
It follows that for $t$ large enough,
$$\displaylines{\qquad
  \P\{\, \exists n\,:\, \epsilon V(t)\leq n\leq V(t)/\epsilon \, ;\,
         S^0_n>t+s_n\,\}
  \hfill\cr\hfill
  {}\leq {V(t)\over \epsilon} 
  \exp\bigl( -r\circ V(t) \theta (1-\epsilon)\bigr)
  \, .
  \qquad\cr}
$$
Still using the large deviations assumption \HypLDP, for $n$ at least 
$V(t)/\epsilon$ and $t$ large enough, we have
$$\eqalign{
  \P\{\, S^0_n>t+s_n\,\}
  &{}\leq \P\{\, S^0_n>s_n\,\} \cr
  &{}\leq e^{-r_n I(1)/2} \, .\cr
  }
$$
Thus, for $t$ large enough, using Lemma \LemmaMasterTh,
$$\eqalign{
  \P\{\, \exists n\, :\, n\geq V(t)/\epsilon\, ;\, S^0_n>t+s_n\,\}
  &{}\leq \sum_{n\geq V(t)/\epsilon} e^{-r_n I(1)/2} \cr
  &{}\lsim \epsilon^{-1} V(t) e^{-r\circ V(t)I(1)/2\epsilon^\rho}
   \, . \cr
  }
$$
Taking $\epsilon$ small enough, it follows that
$$
  \log\P\{\, \exists n \,:\, n\geq V(t)/\epsilon\, ;\, S^0_n>t+s_n \,\}
  \lsim -2r\circ V(t)\theta
$$
as $t$ tends to infinity. Using assumption \HypSmallN, we conclude that
$$
  \log\P\{\, M>t\,\}\lsim -r\circ V(t) \theta
$$
asymptotically.

\noindent{\it Lower bound.} Let $\epsilon$ be a positive real number and let
$x$ be a positive real number such that $x^\rho I(x^{-\gamma}+1)\leq\theta
+\epsilon$. Let $n$ be the integer part of $xV(t)$. Then
$$
  \P\{\, M>t\,\}
  \geq \P\{\, S^0_n>t+s_n\,\} \, .
$$
Using the large deviations hypothesis \HypLDP, we deduce
$$\eqalignno{
  \log\P\{\, M>t\,\}
  &{}\gsim -r_n I\Bigl({t\over s_n}+1\Bigr) \cr
  &{}\sim -r\circ V(t) x^\rho I(x^{-\gamma}+1) \cr
  &{}\gsim -r\circ V(t)(\theta+\epsilon) \, .
  &\equa{EqMasterThB}\cr
  }
$$
Since $\epsilon$ is arbitrary, the first assertion of 
Proposition \MasterTh\ follows.

To prove the second assertion, note that estimate \HypLDP\ with \HypMinUnique\
imply
$$\displaylines{\qquad
  \P\{\, | \calN_t-\tau|>\eta \mid M>t\,\}
  \hfill\cr\hfill
  \leq { \P\{\, \exists n\,:\, |n-\tau V(t)|>\eta V(t) \, ;\, S^0_n>t+s_n\,\}
         \over
         \P\{\, M>t\,\} }
  \qquad\cr
  }
$$
tends to $0$ as $t$ tends to infinity. The second assertion follows.\hfill\qed

\bigskip


\def\prevs{\the\sectionnumber .\the\subsectionnumber .\the\snumber }
\def\preveq{(\the\sectionnumber .\the\subsectionnumber .\the\equanumber)}

\section{Proof of results of Section 2.}
Except if indicated otherwise, we will assume that the mean of the innovations,
$\mu$, is $-1$. Other values of $\mu$ will be dealt with by a scaling argument.

To obtain pleasing expressions, for every positive real number $r$ we
write $\g{r}$ for $\sum_{0\leq i<r} g_i$ and we also write $s_n$ for
the negative of the mean of $S_n$, that is $s_n=\g{n}$  ---~recall
our assumption that $\mu$ is $-1$ until further notice. With the
notation of the previous section, $S_n^0$ is the centered process
$S_n-ES_n=S_n+s_n$. Moreover, $V$ is defined by $s_{\lfloor V(t)\rfloor}\sim t$
as $t$ tends to infinity.

\bigskip

\subsection{Preliminary}%
The following lemma, relating $g_n$ and $\g{n}$ to $g(1-1/n)$ 
will be very useful. It essentially restates Karamata's Tauberian theorem for
power series (Bingham, Goldie and Teugels, 1989, Corollary 1.7.3) and is
proved in Lemma 5.1.1 in Barbe and McCormick (2008). We state it here for 
the sake of making the proof easier to read, for it is fundamental in our 
problem and we will refer to it often.

\Lemma{\label{BS}%
  The following asymptotic equivalences hold as $n$ tends to infinity, 
  uniformly in $x$ in any compact subset of the positive half-line,

  \smallskip

  \noindent (i) $g_{\lfloor nx\rfloor}
  \sim {\ds x^{\gamma-1}\over\ds \Gg}\,{\ds g(1-1/n)\over\ds n}$,

  \smallskip

  \noindent (ii) $\g{nx}
  \sim {\ds x^\gamma\over\ds\Gamma(1+\gamma)}\, g(1-1/n)$.
}

\bigskip

In particular, this implies that $g_n\sim \gamma\g{n}/n$ as $n$ tends to
infinity, so that locally uniformly in any positive $c$,
$$
  g_{\lfloor c V(t)\rfloor} 
  \sim \gamma c^{\gamma-1}{t\over V(t)}
  \eqno{\equa{gEquiv}}
$$
as $t$ tends to infinity.

We introduce the notation $g_{i/n}$ for $\gamma g_i/g_n$ in which the 
subscript $i/n$ has clearly nothing to do with the division of $i$ by $n$ 
but serves as a mnemonic for the division of $g_i$ by $g_n$. In particular,
$g_{n-i/n}$ is $\gamma g_{n-i}/g_n$. Lemma \BS\ asserts that 
$g_{n-i/n}\sim k_\gamma (i/n)$ as $n$ tends to infinity and $i/n$ 
stays bounded away from $1$.

The following easy lemma is recorded for further reference.

\Lemma{\label{UGrowth}
  Let 
  $$
    c_1=\liminf_{n\to\infty} g_n 
    \quad\hbox{and}\quad
    c_2=\limsupn \max_{0\leq i\leq n} g_{i/n} \, .
  $$

  \noindent (i) If $c_1$ is positive, then 
  $U\lsim \Id/c_1\Gamma(\gamma)$ at infinity.

  \smallskip

  \noindent (ii) Assume that $c_2$ is finite. If the sequence 
  $(g_n)_{n\geq 0}$ is bounded, then $U\lsim c_2\Id/\gamma\max_{i\geq 0} g_i$; 
  otherwise $U=o(\Id)$ at infinity.
}

\bigskip

\Proof To prove (i), Lemma \BS.i ensures that
$$
  g(1-1/n) \gsim \Gamma(\gamma)c_1 n 
$$
as $n$ tends to infinity. Therefore, since $g$ is regularly varying,
$$
  \Id\sim g(1-1/U)\gsim c_1\Gamma(\gamma) U
$$
at infinity, and the result follow.

To prove (ii), let $c$ be a number greater than $c_2$ and let $k$ be an integer
such that $g_k$ is positive. Then $g_n\geq \gamma g_k/c$
for any $n$ larger than some $n_0$. Therefore, on $[\,0,1)$, the function $g$
is bounded from below by a polynomial of degree $n_0$ plus the function
$$
  {\gamma g_k\over c} \sum_{n\geq n_0} x^n  
  = {\gamma g_k\over c(1-x)} x^{n_0} \, .
$$
This implies that
$$
  \Id\sim g(1-1/U)\gsim (\gamma g_k/c) U
$$
at infinity. Since $c$ and $k$ are arbitrary, this prove assertion 
(ii).\hfill\qed

\bigskip

Our next lemma is perhaps the heart of the proof, which ultimately relies on
approximation of Riemann sums by Riemann integral, a modicum of regular 
variation, and the exponential form of Markov's inequality.

We define the sequence of probability 
measures
$$
  {\mit\Gamma}_n = n^{-1}\sum_{1\leq i\leq n} \delta_{(i/n,g_{n-i/n})} \, , 
  \quad 
  n\geq 1 \, .
$$

\Lemma{\label{GammaNConverges}
  The sequence of probability measures $({\mit\Gamma}_n)_{n\geq 1}$ converges
  weakly$*$ to the measure $\int_0^1\delta_{(u,k_\gamma(u))}\d u$.
}

\bigskip

\Proof Let $f$ be a nonnegative continuous and bounded function on 
$[\,0,1\,]\times \RR$. We write $|f|_{[0,1]\times\RR}$ for its supremum on
the strip $[\,0,1\,]\times\RR$. Let $\epsilon$ be a positive real number 
less than $1$. Note that
$$
  n^{-1} \sum_{(1-\epsilon)n< i\leq n} f(i/n,g_{n-i/n})
  \leq \epsilon |f|_{[0,1]\times \RR} \, .
$$
Uniformly in $i$ between $1$ and $(1-\epsilon)n$, Lemma \BS\ shows
that $g_{n-i/n}\sim k_\gamma(i/n)$. Thus, since the measure 
$n^{-1}\sum_{1\leq i\leq n}\delta_{i/n}$ converges weakly$*$ to the Lebesgue 
measure on $[\,0,1\,]$ as $n$ tends to infinity, 
$$
  \limn n^{-1}\sum_{1\leq i\leq (1-\epsilon)n} f(i/n,g_{n-i/n})
  = \int_0^{1-\epsilon} f\bigl( u,k_\gamma(u)\bigr) \d u \, ,
$$
and the result follows.\hfill\qed

\bigskip

In order to simplify the notation during the proof and make a later scaling
argument easier to follow, we write $\varphi_0$ for the moment generating
function of the centered random variable $(X_1/{-\mu})+1$. Furthermore, we
write
$$
  J_0(a)=\sup_{\lambda>0} 
  \Bigl(a\lambda-\int_0^1\log\varphi_0\bigl(\lambda k_\gamma(u)\bigr)\d u 
  \Bigr)
  \, . 
  \eqno{\equa{JZeroDef}}
$$
The equality $\varphi_0(\lambda)=e^\lambda\varphi(\lambda/{-\mu})$, valid
for all $\lambda$ positive, yields
$$\eqalignno{
  J_0(a+1)
  &{}=\sup_{\lambda>0}\Bigl( (a+1)\lambda 
   -\int_0^1\log\varphi_0\bigr(\lambda k_\gamma (u)\bigl)\d u \Bigr) \cr
  &{}=J(-\mu a) \, . 
  &\equa{JJeq}\cr
  }
$$

\bigskip

\subsection{Proof of Theorem \ThCaseA}%
Recall that except if specified otherwise, we consider the mean $\mu$ to
be $-1$. Also, throughout this subsection, we assume that the the hypotheses
of Theorem \ThCaseA\ hold, even if this is not specified.

The proof is based on a large deviations estimate which is the analogue for
$(g,F)$-process of the classical estimate of Chernoff for the sample mean.
The proof requires several lemmas.

Our first lemma will be useful in taking limits in various sums involving
the moment generating function.

\Lemma{\label{LDPLemma}
  Let $h$ be a continuous function on the nonnegative half-line. 
  Assume $g$ satisfies the assumption of Theorem \ThCaseA.
  If $\limn \max_{0\leq i\leq n} g_{i/n}=\infty$,
  assume further that $h$ is regularly varying of index $\beta$ less 
  than $1/(1-\gamma)$. Then, locally uniformly in $\lambda$ in $(0,\infty)$,
  $$
    \limn n^{-1}\sum_{1\leq i\leq n} h(\lambda g_{n-i/n})
    =\int_0^1h\bigl(\lambda k_\gamma(u)\bigr)\d u \, ,
  $$
  and this limit is finite.
}

\bigskip

\noindent When $\limn\max_{0\leq i\leq n}g_{i/n}=\infty$ and $\gamma$ is $1$, 
the condition on $h$ should simply be read as $h$ is regularly varying of some
positive index.

\bigskip

\Proof Note first that in both cases, 
$\limeps \int_0^\epsilon h(\lambda u^{\gamma-1})\d u=0$ and the 
integral involved in the limit in the lemma is indeed finite.

Once the limit is established for a fixed $\lambda$, it will be clear
that using the uniform convergence theorem for regularly varying
functions, the limit is locally uniform in $\lambda$. Thus, up to changing the
function $h$, it suffices to prove the result only when $\lambda$ is $1$.

For any positive real number $c$ we define the function 
$h_c=h(\,\cdot\,\wedge c)$. These functions are continuous and bounded.

If $\limsupn \max_{0\leq i\leq n}g_{i/n}$ is finite, we take $c$ to be twice
this limit, so that for $n$ large enough, $\max_{0\leq i\leq n} g_{i/n}$ is
at most $c$. Then
$$
  n^{-1}\sum_{1\leq i\leq n} h(g_{n-i/n})
  =\int h_c(x)\d{\mit\Gamma}_n(v,x)
$$
and the result follows from Lemma \GammaNConverges\ and the local
uniform continuity in $\lambda$ of the functions $h_c(\lambda\,\cdot\,)$.

If $\limn\max_{0\leq i\leq n}g_{i/n}$ is infinite, let $\epsilon$ be a positive
real number less than $1$. Since $\limn\max_{\epsilon n\leq i\leq n} g_{i/n}
=\epsilon^{\gamma-1}$ is finite, it suffices to prove that
$$
  \lim_{\epsilon\to 0} \limsupn n^{-1}\sum_{0\leq i\leq \epsilon n} h(g_{i/n})
  = 0 \, .
  \eqno{\equa{eqLDPLemma}}
$$

Let $\delta$ be a positive real number, small enough so that
\hbox{$\beta(1-\gamma+\delta)$} is less than $1$. Lemma \BS\ and Potter's
bound (Bingham, Goldie and Teugels, 1989, Theorem 1.5.6) show that
there exists $n_0$ such that whenever $n$ and $i$ are at least $n_0$,
then $g_{i/n}$ is at most $2\gamma(i/n)^{\gamma-1-\delta}$. Since $h$ is
regularly varying with positive index $\beta$, it is asymptotically
equivalent to a nondecreasing function (\hbox{Bingham}, Goldie and Teugels,
1989, Theorem 1.5.3). Hence, provided $n$ is large enough,
$$\eqalign{
  n^{-1}\sum_{n_0\leq i<n\epsilon} h(g_{i/n})
  &{}\leq 2n^{-1}\sum_{n_0<i\leq n\epsilon} 
    h\bigl( 2\gamma (i/n)^{\gamma-1-\delta}\bigr)\cr
  &{}\leq 4\int_0^\epsilon h(2\gamma u^{\gamma-1-\delta}) \d u \, .\cr
  }
$$
Moreover, for $M$ large enough and as $n$ tends to infinity,
$$\eqalign{
  n^{-1}\sum_{0\leq i<n_0}h(g_{i/n})
  &{}\leq n^{-1} n_0 h(M/g_n) \cr
  &{}=O\Bigl(n^{-1} h\bigl(n/g(1-1/n)\bigr)\Bigr)\, . \cr
  }
$$
This bounds tends to $0$ as $n$ tends to infinity since the function $x\mapsto
x^{-1}h\bigl(x/g(1-1/x)\bigr)$ is regularly varying with negative 
index $\beta(1-\gamma)-1$. This proves \eqLDPLemma.\hfill\qed

\bigskip

Equipped with Lemma \LDPLemma, we can prove the following large deviation
principle. Recall that we assume for the time being that the distribution
function $F$ has mean $-1$. We write $F_0$ for the cumulative distribution 
function $F(\cdot -1)$. As the subscript indicates, its mean is $0$.

\Proposition{\label{LDP}
  Let $(S^0_n)_{n\geq 0}$ be a centered $(g,F_0)$-process.
  Under the assumptions of Theorem \ThCaseA,
  for any nonnegative $x$,
  $$
    \limn n^{-1}\log \P\{\, S^0_n>  x\g n\,\}=-J_0(x) \, .
  $$
  Moreover, the limit is locally uniform in $x$ on the set where $J_0$ is 
  finite.}%

\bigskip

\Proof The proof is modeled after the standard one for the mean. We
will concentrate on proving a pointwise version in $x$ because the
following purely analytical argument gives the local uniformity. If
the pointwise result holds, it asserts that the sequence of
nonincreasing functions $(n^{-1}\log\P\{\,S_n> \g n \,\cdot\,\,\})_{n\geq 1}$ 
converges to the function $-J_0$; since the limit is
continuous, and monotone as a limit of monotone functions, the
convergence is locally uniform (see Rudin, 1976, chapter 7, exercise
13).

\noindent{\it Upper bound.} Let $\lambda$ be a positive real number. Using
the exponential Markov inequality,
$$
  \P\Bigl\{\, {\gamma S^0_n\over g_n}> x{\gamma\g n\over g_n} \,\Bigr\}
  \leq \exp\biggl( -\lambda x{\gamma\g n\over g_n} +\sum_{1\leq i\leq n}
  \log\varphi_0(\lambda g_{n-i/n})\biggr) \, .
$$
This implies that
$$\displaylines{\qquad
  \limsupn n^{-1}\log \P\{\, S^0_n> x\g n\,\}
  \hfill\cr\hfill
  {}\leq \limsupn \biggl( -\lambda x{\gamma\g n\over ng_n} 
  +n^{-1}\sum_{1\leq i\leq n} \log\varphi_0(\lambda g_{n-i/n}) \biggr) \, .
  \qquad\equa{LDPA}
  \cr}
$$
Note that Lemma \BS\ implies that
$$
  {\gamma\g n\over ng_n}
  \sim 1
$$
while Lemma \LDPLemma\ yields, in view of the fact that $\log\varphi_0$
is regularly varying with index $\beta$, the conjugate exponent to $\alpha$,
and $\alpha\gamma>1$, that
$$
  n^{-1}\sum_{1\leq i\leq n} \log\varphi_0 (\lambda g_{n-i/n})
  \sim \int_0^1\log\varphi_0\bigl(\lambda k_\gamma(u)\bigr)\d u
$$
as $n$ tends to infinity. Therefore, the right hand side of \LDPA\ tends
to
$$
  -\lambda x +\int_0^1\log\varphi_0 \bigl(\lambda k_\gamma (u)\bigr)\d u
$$
as $n$ tends to infinity. The infimum of this upper bound over 
all $\lambda$ positive is $-J_0(x)$.

\noindent{\it Lower bound.} We write $\II\{~\}$ for the indicator function
of a set. For any fixed $\lambda$, the equality
$$\eqalignno{
  &\hskip-20pt \P\Bigl\{\, {S_n^0\over g_n}> x{\g n\over g_n}\,\Bigr\}\cr
  &{}= \Bigl(\prod_{1\leq i\leq n}\varphi_0(\lambda g_{n-i/n}) 
    e^{-\lambda x\gamma\g n/g_n}\Bigl)
    \int \II\Bigl\{\, {S_n^0\over g_n}\geq 
                      x{\g n\over g_n}\,\Bigr\}\cr
  &\qquad\qquad{}\times e^{-\lambda\sum_{1\leq i\leq n} g_{n-i/n}x_i} 
               e^{\lambda x\gamma\g n/g_n} \cr
  &\qquad\qquad{}\times { e^{\lambda\sum_{1\leq i\leq n}g_{n-i/n}x_i}
                   \over\prod_{1\leq i\leq n}
                   \varphi_0 (\lambda g_{n-i/n}) }  
          \d F_0(x_1)\ldots \d F_0(x_n)
  &\equa{LDPUBa}\cr}
$$
holds. We write $Q_{\lambda,n}$ for the image measure of the probability 
measure
$$
  { e^{\lambda\sum_{1\leq i\leq n}g_{n-i/n}x_i}
    \over 
    \prod_{1\leq i\leq n}\varphi_0(\lambda g_{n-i/n}) }
  \d F_0(x_1)\ldots \d F_0(x_n)
$$
through the map
$$
  (x_1,\ldots , x_n)\mapsto 
  \sum_{1\leq i\leq n} g_{n-i/n}x_i-x\gamma{\g n\over g_n}  \, . 
$$
Writing $R$ for
a random variable having distribution $Q_{\lambda,n}$, equality \LDPUBa\ is
equivalent to
$$\displaylines{\quad
  \P\{\, S^0_n> x\g n \,\} 
  \hfill\cr\hfill
  = \prod_{1\leq i\leq n}\varphi_0 (\lambda g_{n-i/n})\,\, 
   e^{-\lambda x\gamma\g n /g_n}
  \E\II\{\, R\geq 0\,\} e^{-\lambda R}\, . 
  \quad\equa{LDPUBaa}\cr
  }
$$
The moment generating function of $R$ evaluated at $\zeta$ is
$$\eqalignno{
  \E e^{\zeta R}
  &{}={ \E e^{(\lambda+\zeta)\gamma S^0_n/g_n - \zeta x\gamma \g n/g_n}
        \over
        \prod_{1\leq i\leq n} \varphi_0 (\lambda g_{n-i/n}) } \cr
  &{}=e^{-\zeta  x\gamma \g n/g_n} \prod_{1\leq i\leq n} 
      { \varphi_0 \bigl( (\lambda+\zeta)g_{n-i/n}\bigr) 
        \over 
        \varphi_0 (\lambda g_{n-i/n}) } \, .
  &\equa{LDPUBb}\cr
  }
$$
In particular, taking its logarithmic derivative at $0$, the expectation
of $R$ is
$$
  \E R=-{x\gamma\g n \over g_n} 
  +\sum_{1\leq i\leq n} g_{n-i/n}m_0(\lambda g_{n-i/n})
  \, ,
$$
while its variance is
$$
  \sigma_n^2=\sum_{1\leq i\leq n} g_{n-i/n}^2m_0'(\lambda g_{n-i/n}) \, .
$$
Using Lemma \LDPLemma, we obtain that, locally uniformly in $\lambda$,
$$
  \E R\sim n\Bigl( -x 
                   +\int_0^1 k_\gamma(u) m_0\bigl(\lambda k_\gamma(u)\bigr)\d u
            \Bigr)
  \, .
$$
Since $m_0'$ is assumed regularly varying when $\limsup_{n\to\infty}g_{i/n}$
is infinite,
$$
  \sigma_n^2\sim n\int_0^1 k_\gamma^2(u)
  m_0'\bigl(\lambda k_\gamma(u)\bigr)\d u \, ,
$$
and the integral involved in this asymptotic equivalence is well
defined ---~see the discussion following Theorem \ThCaseA, where we
showed that in fact $m_0'$ is regularly varying with index $\beta-2$
with $\beta(1-\gamma)<1$.

We consider $\lambda$, depending on $n$, therefore written $\lambda_n$ from
now on, such that the expected value of $R$ vanishes. Such sequence
exists since the standard assumption guaranties that $m$ is onto 
the nonnegative real line. Since $m_0$ is monotone, this sequence $\lambda_n$ 
converges to the root of
$$
  -x +\int_0^1 k_\gamma(u)
  m_0\bigl(\lambda k_\gamma(u)\bigr)\d u = 0 \, .
$$
Then \LDPUBaa\ implies
$$\displaylines{
  \P\{\, S^0_n> x\g n\,\}
  \hfill\cr
    {}=\prod_{1\leq i\leq n}\hskip -3pt \varphi_0(\lambda_n g_{n-i/n}) 
         e^{-\lambda_n x\gamma\g n/g_n}
         \E \II\{\, R\geq 0\,\} e^{-\lambda_n \sigma_nR/\sigma_n} 
  \hfill\cr
    {}\geq\prod_{1\leq i\leq n}\hskip -3pt\varphi_0(\lambda_n g_{n-i/n})
         e^{-\lambda_n x\gamma\g n/g_n}
           e^{-\lambda_n \sigma_n M} Q_{\lambda_n,n}[\, 0,M\sigma_n\,] \, .
    \hfill\equa{LDPUBc}
  \cr}
$$
Expression \LDPUBb\ shows that the logarithm of the moment generating 
function of $R/\sigma_n$ at $\zeta$ is
$$
  -{\zeta x\gamma\g n\over \sigma_n g_n} 
  + \sum_{1\leq i\leq n} \log\varphi_0 \bigl((\lambda_n
  +\zeta/\sigma_n)g_{n-i/n})\bigr)-\log\varphi_0 (\lambda_n g_{n-i/n}) \, .
$$
Using the mean value theorem and given our choice of $\lambda_n$, there exists
some $\eta_{i,n}$ between $0$ and $1$ such that this 
logarithm is
$$
  {\zeta^2\over 2\sigma_n^2} \sum_{1\leq i\leq n} g_{n-i/n}^2 m_0'
  \bigl( (\lambda_n+\eta_{i,n}\zeta/\sigma_n)g_{n-i/n}\bigr)\, .
  \eqno{\equa{LDPUBd}}
$$
The same argument as in Lemma \LDPLemma\ shows that \LDPUBd\ tends 
to $\zeta^2/2$ as $n$ tends to infinity. Therefore, $R/\sigma_n$ has a
standard Gaussian limiting distribution as $n$ tends to infinity, and 
the right hand side of \LDPUBc\
is asymptotically equivalent to
$$\displaylines{\qquad
  \exp\Bigl( \sum_{1\leq i\leq n}\log
  \varphi_0 (\lambda_n g_{n-i/n})-\lambda_n x\gamma\g n/g_n
  \Bigr) e^{O(\sqrt n)}
  \hfill\cr\hfill
  {}=\exp\Bigl( -nJ_0(x)\bigl( 1+o(1)\bigr)\Bigr)
  \qquad\cr}
$$
as $n$ tends to infinity. The result follows.\hfill\qed

\bigskip

To prove Theorem \ThCaseA\ requires a couple of more
lemmas related to the function $J_0$.

\Lemma{\label{IZeroPositive}
  The function $J_0$ is positive on the positive half-line.
  }

\bigskip

\Proof  Let $J_0^*$ be the function
$$
  J_0^*(\lambda)=\int_0^1\log\varphi_0 \bigl(\lambda k_\gamma(u)\bigr)\d u \, .
$$
Both $J_0^*$ and ${J_0^*}'$ vanish at the origin, while
$$
  {J_0^*}''(0)=\Var X_1\int_0^1 k_\gamma^2(u) \d u
$$
is positive. In particular, taking $\lambda$
to be $x/{J_0^*}''(0)$, we see that as $x$ tends to $0$,
$$\eqalign{
  J_0(x)
  &{}\geq \lambda x -{\lambda^2\over 2} {J_0^*}''(0) +o(\lambda^2)\cr
  &{}\geq {x^2\over 2{J_0^*}''(0)} +o(x^2) \, . \cr
  }
$$
Thus, $J_0$ is positive on an open interval with left endpoint the origin. 
Since $J_0$ is a supremum of nondecreasing functions of $x$ it is also 
nondecreasing and the result follows.\hfill\qed

\Lemma{\label{IMin}
  For any positive real number $c$, the function 
  $x\in [\,0,\infty)\mapsto xJ_0(x^{-\gamma}+c)$ tends to infinity
  at $0$ and infinity. Moreover, it reaches its minimum at a positive argument.
}

\bigskip

\Proof Let $c$ be a positive real number. Lemma \IZeroPositive\ ensures 
that $J_0(c)$ is positive. It follows that $xJ_0(x^{-\gamma}+c)$ tends to 
infinity with $x$.

Assume that $\gamma$ is at least $1$. Since for any 
positive $\theta$ the
inequality $J_0(x)\geq x\theta-J_0^*(\theta)$ holds, we see 
that $J_0$ ultimately grows faster than any multiple of the identity. 
Thus, $xJ_0(x^{-\gamma}+c)$ tends to 
infinity as $x$ tends to $0$, and this proves the lemma in this case.

Assume that $\gamma$ is less than $1$. The assumption $\alpha\gamma>1$
ensures that $-\log\oF$ is
regularly varying of index $\alpha$ greater than $1$. By Kasahara's
(1978) Tauberian theorem (Bingham, Goldie and Teugels, 1989, Theorem
4.12.7), $\log\varphi_0$ is regularly varying of index $\beta$, the
conjugate exponent to $\alpha$.  This implies that $J_0^*$ is also
regularly varying of index $\beta$ at infinity. By Bingham and
Teugels's (1975) theorem (see Bingham, Goldie and Teugels, 1989,
Theorem 1.8.10), this implies that $J_0$ is regularly varying of index
$\alpha$. Since $\alpha\gamma$ is greater than $1$, it then follows
that $x J_0(x^{-\gamma}+c)$ tends to infinity as $x$ tends to $0$
(Bingham, Goldie and Teugels, 1989, Proposition 1.3.6). This proves
the first part of the lemma.

The second part of the lemma follows, because the 
function $xJ_0(x^{-\gamma}+c)$ is
continuous on the positive half-line.\hfill\qed

\bigskip

Our next lemma shows that the process $S_n^0$ is unlikely to reach a high 
threshold $t$ before a time of order $V(t)$.

\bigskip

\Lemma{\label{ThASmallIndices}
  The following holds,
  $$
    \lim_{\epsilon\to 0}\limsupt {1\over V(t)} \log 
    \P\{\, \exists n \, :\, n\leq\epsilon V(t)\, ;\, S^0_n>t\,\} 
     =-\infty \, .
  $$
}

\Proof We distinguish according to whether $\max_{0\leq i\leq n} g_{i/n}$ 
remains bounded or not.

Assume first that $\limsupn\max_{0\leq i\leq n} g_{i/n}$ is some
finite positive number $c$. Necessarily, $\gamma$ is at least $1$. In
that case, \gEquiv\ implies
$$
  \max_{0\leq i\leq \epsilon V(t)} g_i
  \lsim {c\over\gamma} g_{\lfloor\epsilon V(t)\rfloor}
  \sim c\epsilon^{\gamma-1}{t\over V(t)}
$$
as $t$ tends to infinity. In particular, uniformly in $i$ nonnegative and
at most $\epsilon V(t)$, and as $t$ tends to infinity, 
$t/g_i\gsim V(t)\epsilon^{1-\gamma}/c$. Moreover, Lemma \LDPLemma\ shows
that
$$
  \sum_{1\leq i\leq n} \log\varphi_0(\lambda g_{n-i/n})
  \lsim n\int_0^1\log\varphi_0\bigl(\lambda k_\gamma(u)\bigr) \d u
$$
as $n$ tends to infinity. Then, using the Markov exponential inequality,
for any fixed positive $\lambda$, for any $n$ large enough and at most
$\epsilon V(t)$,
$$\eqalign{
  \log\P\{\, S^0_n>t\,\}
  &{}\leq -\lambda {\gamma t\over g_n} + \sum_{1\leq i\leq n} 
   \log\varphi_0(\lambda g_{n-i/n})\cr
  &{}\leq -\lambda \gamma {V(t)\over 2c} \epsilon^{1-\gamma} 
   + 2n\int_0^1\log\varphi_0 \bigl(\lambda k_\gamma(u)\bigr)\d u \cr
  }
$$
provided $t$ is large enough, $n$ is large enough and less 
than $\epsilon V(t)$.

Since $\gamma$ is at least $1$, for $n$ at most $\epsilon V(t)$, this upper
bound is at most
$$ 
  -V(t) \Bigl( \lambda \gamma {\epsilon^{1-\gamma}\over 2c}
               -2\epsilon\int_0^1\log\varphi_0(\lambda k_\gamma(u)\d u\Bigr) 
  \, .
$$
It can be made smaller than any negative multiple of $V(t)$
by first taking $\lambda$ positive and then $\epsilon$ small enough. Hence,
there exists $n_0$ such that
$$
  \limeps\limsupt \max_{n_0\leq n\leq\epsilon V(t)}{1\over V(t)}
  \log\P\{\, S^0_n>t\,\} = -\infty \, .
  \eqno{\equa{UBdd}}
$$

For $n$ at most $n_0$, recalling that the mean of $X_i$ is $-1$, we have,
since $t$ is positive,
$$\eqalign{
  P\{\, S^0_n>t\,\}
  &{}\leq  P\{\, n_0\max_{0\leq i\leq n_0} g_i\max_{1\leq i\leq n_0} (X_i+1)>t
              \,\} \cr
  &{}\leq n_0\oF_0\Bigl( {t\over n_0\max_{0\leq i\leq n_0} g_i} \Bigr) \, . \cr
  }
$$
Since the moment generating function of $F_0$ is finite on the nonnegative
half-line, Chernoff's inequality implies that $-\log\oF\gg \Id$ at infinity.
Lemma \UGrowth\ shows that in the present case, $U$ grows at most like
a multiple of the identity at infinity. This implies that the 
function $U^{-1}\log\oF$ tends to 
minus infinity at infinity. We conclude that \UBdd\ holds with $n_0$ being $1$.

We now consider the case 
where $\max_{0\leq i\leq n} g_{i/n}$ tends to infinity with $n$. In this case,
the sequence $(g_n)_{n\geq 0}$ converges to $0$, and, 
for any $\eta$ positive, $\log\varphi_0\lsim \Id^{\beta+\eta}$ at infinity.
Again, we use the exponential Markov inequality
$$
  \log\P\{\, S^0_n>t\,\}
  \leq -\lambda t+\sum_{0\leq i< n} \log\varphi_0(\lambda g_i) \, ,
  \eqno{\equa{UBddd}}
$$
taking now $\lambda$ of the form $cV(t)/t$ for some positive constant $c$
to be determined.

Since the standard assumption ensures that $(g_n)_{n\geq 0}$ is asymptotically
equivalent to a monotone sequence, 
$\min_{0\leq i\leq \epsilon V(t)} g_i\gsim g_{\lfloor \epsilon V(t)\rfloor}$
as $t$ tends to infinity. Using \gEquiv, it follows that
$\lambda g_i\gsim c\gamma\epsilon^{\gamma-1}$ is large whenever $c$ is large 
and $\epsilon$ is small. Thus, provided $c$ is large enough, $\epsilon$ is 
small enough and $n$ is at most $\epsilon V(t)$,
$$\eqalignno{
  \sum_{0\leq i<n}\log\varphi_0(\lambda g_i)
  &{}\leq 2\sum_{0\leq i<n} (\lambda g_i)^{\beta+\eta} \cr
  &{}\leq 2 \Bigl( c{V(t)\over t}\Bigr)^{\beta+\eta} 
   \sum_{0\leq i\leq\epsilon V(t)} g_i^{\beta+\eta} \, .
  &\equa{UBdddd}\cr
  }
$$
Since $\beta(\gamma-1)+1$ is positive, 
$$
  \sum_{0\leq i<n} g_i^{\beta+\eta}
  \sim {n\over1+(\gamma-1)(\beta+\eta)}
       \Bigl( {g(1-1/n)\over\Gg n}\Bigr)^{\beta+\eta}
$$
as $n$ tends to infinity, and the bound \UBdddd\ is at most
$$
  2 ( \gamma c)^{\beta+\eta} \epsilon^{(\gamma-1)(\beta+\eta)+1}
  {V(t)\over 1+(\gamma-1)(\beta+\eta)} \, .
$$
For any fixed $\epsilon$, the upper bound \UBdddd\ 
can be made less than any a priori given negative number times $V(t)$
by taking $c$ large enough. This 
proves the lemma.\hfill\qed

\bigskip

\noindent{\bf Proof of Theorem \ThCaseA.}
Comparing \HypLDP\ with Proposition \LDP, we may take $r$ to be the identity,
so that $\rho$ is $1$; furthermore, still referring to assumption \HypLDP\ and 
Proposition \LDP, we see that $I(x)=J_0(x)$. Using \JJeq, and since the
mean $\mu$ is $-1$, it follows that
$\theta$, as defined in \ThetaDef, is
$$
  \theta
  =\inf_{x>0} x J_0(x^{-\gamma}+1)
  =\inf_{x>0}xJ(x^{-\gamma}) \, .
  \eqno{\equa{ThetaJZero}}
$$
Assumptions needed to apply the first part of Proposition \MasterTh\
are satisfied thanks to Lemma \BS, Proposition \LDP\ and Lemma
\ThASmallIndices. Thus Proposition \MasterTh\ yields Theorem \ThCaseA\ when
$\mu$ is $-1$.

To obtain the result when $\mu$ is different than $-1$, we index in an obvious
way all quantities by the mean $\mu$ in parentheses. We then have, assuming
now that $X_i$ has arbitrary mean $\mu$,
$$
  \varphi_{(\mu)}(\lambda)
  = \varphi_{(-1)}(-\mu\lambda) \, .
  \eqno{\equa{PhiScale}}
$$
This implies, for any positive $a$,
$$
  J_{(\mu)}(a)=J_{(-1)}(a/{-\mu}) \, .
  \eqno{\equa{JScale}}
$$
Writing $X_{(\mu),i}$ for $X_i$ when the mean is $\mu$, we take $X_{(\mu),i}
=-\mu X_{(-1),i}$, which yields $M_{(\mu)}=(-\mu)M_{(-1)}$. Thus,
$$\eqalignno{
  {1\over V(t)} \log\P\{\, M_{(\mu)}>t\,\}
  &{}={V(t/{-\mu})\over V(t)}
     {1\over V(t/{-\mu})} \log\P\{\, M_{(-1)}>t/{-\mu}\,\} \cr
  &{}\sim -(-\mu)^{-1/\gamma} \inf_{x>0} x J_{(-1)}(x^{-\gamma}) \, . \cr
  &{}= -\inf_{x>0} x J_{(\mu)}(x^{-\gamma}) \, . 
  &\qed\cr
  }
$$

\bigskip

\subsection{Proof of Theorem \ThCaseC}
As for Theorem \ThCaseA, we first prove Theorem \ThCaseC\ when $\mu$ is $-1$,
which we assume from now on.

Our first lemma is an analogue of Lemma \LDPLemma\ but in the context 
of Theorem \ThCaseC.

Note that in the context of Theorem \ThCaseC, the
conditions $\alpha\gamma<1$ and $\gamma>1/2$ force $\alpha$ to be less
than $2$. Therefore, its conjugate exponent, $\beta$, is greater than $2$.

\Lemma{\label{LDPBLemma}
  Let $\lambda$ be a regularly varying function of index greater than 
  $(2\gamma-1)/(\beta-2)$ and set $\lambda_n=\lambda(n)$. 
  Under the assumptions of Theorem \ThCaseC,
  $$
    \sum_{0\leq i<n}\log\varphi_0(\lambda_ng_i)
    \sim\log\varphi_0(\lambda_n)|g|_\beta^\beta
  $$
  as $n$ tends to infinity.
}

\bigskip

\Proof
Let $\sigma^2$ be the variance of $X_i$. Since $\log\varphi_0\sim
\Id^2\sigma^2/2$ at the origin, for any positive $R$
there exists a positive $c$ such that $\log\varphi_0\leq c\Id^2$ 
on $[\,0,R\,]$. Using Lemma \BS, this implies
$$\eqalignno{
  \sum_{0\leq i<n}\II\{\, \lambda_n g_i\leq R\,\}\log\varphi_0(\lambda_n g_i)
  &{}\leq c\lambda_n^2\sum_{0\leq i<n} g_i^2 
   \cr
  &{}\asymp \lambda_n^2 {g(1-1/n)^2\over n} \, .
  &\equa{LDPBLemmaEqA}\cr
  }
$$
Let $\epsilon$ be a positive real number. Using Potter's bound, we see that
provided $\lambda_n$ and $\lambda_n g_i$ are large enough, and provided
that $i$ is large enough for $g_i$ to be less than $1$,
$$
  {1\over 2} g_i^{\beta+\epsilon}
  \leq { \log\varphi_0(\lambda_n g_i)\over\log\varphi_0 (\lambda_n) }
  \leq 2 g_i^{\beta-\epsilon} \, .
$$
By a standard regular variation theoretic argument, this implies
$$\eqalignno{
  \sum_{0\leq i<n}\II\{\, \lambda_n g_i>R\,\}\log\varphi_0(\lambda_n g_i)
  &{}\sim \log\varphi_0(\lambda_n) \sum_{0\leq i<n} g_i^\beta
   \II\{\, \lambda_n g_i>R\,\} \cr
  &{}\sim \log\varphi_0 (\lambda_n) |g|_\beta^\beta
  &\equa{LDPBLemmaEqB}\cr
  }
$$
as $n$ tends to infinity ---~recall that $|g|_\beta$ is finite here, since
$\beta(1-\gamma)>1$.

Write $\rho$ for the index of regular variation of $\lambda$.
Since $\log\varphi_0\circ\lambda$ is regularly varying of index $\beta\rho$,
and $\lambda^2 g(1-1/\Id)^2/\Id$ is regularly varying of index 
$2\rho+2\gamma-1$, our assumption that $\rho$ is greater than  
$(2\gamma-1)/(\beta-2)$
ensures that the right hand side of \LDPBLemmaEqB\ dominates the right hand
side of \LDPBLemmaEqA, and the result holds.\hfill\qed

\bigskip

We define the Cram\'er transform of the centered random variables,
$$
  I_0(x)=\sup_{\lambda>0} \lambda x-\log\varphi_0(\lambda) \, .
$$
Recall that we assume that $\mu$ is $-1$. We then write $F_0$ for 
the distribution of the centered random variable $X_i+1$.

We can now state and prove the following large deviations inequality.

\Proposition{\label{LDPB}
  Let $(S^0_n)_{n\geq 0}$ be a centered $(g,F_0)$-process.
  Under the assumptions of Theorem \ThCaseC, for any positive $x$,
  $$
    \log\P\{\, S^0_n>x\g n\,\}
    \lsim  -{x^\alpha\over |g|_\beta^\alpha} I_0(\g n)
  $$
  as $n$ tends to infinity.
}

\bigskip

\Proof Recall that under the assumptions
of Theorem \ThCaseC, $\log\varphi_0$ is regularly varying of index $\beta$
and $I_0$ is regularly varying of index $\alpha$ ---~see the proof of 
Lemma \IMin. Define
$$
  \lambda(t)
  ={x^{1/(\beta-1)}\over |g|_\beta^\alpha} m_0^\leftarrow (\g t) \, .
$$
This function is regularly varying of positive index $\gamma/(\beta-1)$. 
We define $\lambda_n$ as $\lambda(n)$.
Using the exponential form of Markov's inequality and Lemma \LDPBLemma 
---~applicable since the inequality $\alpha\gamma<1$ implies 
$\gamma/(\beta-1)>(2\gamma-1)/(\beta-2)$~---
$$\eqalign{
  \log\P\{\, S^0_n>x\g n\,\}
  &{}\leq -\lambda_n x \g n 
            + \sum_{0\leq i<n}\log\varphi_0(\lambda_n g_i) \cr
  &{}\leq -{x^\alpha\over |g|_\beta^\alpha} \g n 
            m_0^\leftarrow\circ \g n \cr
  &\quad\quad{}+|g|_\beta^\beta {x^\alpha\over |g|_\beta^{\alpha\beta}}
           \log\varphi_0\circ m_0^\leftarrow(\g n)\bigl(1+o(1)\bigr) 
           \, . \cr
  }
$$

Since by regular variation $\log\varphi_0\sim \Id m_0/\beta$ at infinity, 
the right hand side of the above upper bound is asymptotically equivalent to
$$
  {x^\alpha\over |g|_\beta^\alpha} (\Id m_0^\leftarrow)(\g n)
  (-1+1/\beta) \, .
  \eqno{\equa{LDPBEqA}}
$$
Upon noting that the maximizing value of $\lambda$ in the definition 
of $I_0(\cdot)$ is $m_0^\leftarrow(\cdot)$, the chain rule yields
$I_0'=m_0^\leftarrow$. Therefore, $\Id m_0^\leftarrow\sim \alpha I_0$
at infinity. We obtain that \LDPBEqA\ is asymptotically equivalent to
$-x^\alpha |g|_\beta^{-\alpha} I_0(\g n)$ as $n$ tends
to infinity. This proves Proposition \LDPB.\hfill\qed

\bigskip

Our next result is yet another large deviations inequality. Its statement is
suitable for our application, though its proof gives a somewhat more precise
estimate.

\Proposition{\label{LDPC}
  For any positive real number $\zeta$,
  $$
    \max_{1\leq n<\zeta V(t)} \log\P\{\, S^0_n>t\,\} 
    \lsim -{I_0(t)\over |g|_\beta^\alpha}
  $$
  as $t$ tends to infinity.
  }

\bigskip

\Proof Let $\lambda(t)= m_0^\leftarrow (t) |g|_\beta^{-\alpha}$. 
The exponential Markov inequality implies
$$
  \log\P\{\, S^0_n>t\,\}
  \leq -\lambda(t) t+\sum_{0\leq i<n}
    \log\varphi_0\bigl(\lambda(t) g_i\bigr) \, . 
  \eqno{\equa{LDPCEqA}}
$$
Using Potter's bound and regular variation of $\log\varphi_0$, there exists
a positive $R$ such that uniformly in $n$ positive and less than $\zeta
V(t)$,
$$
  \sum_{0\leq i<n}\log\varphi_0\bigl(\lambda(t)g_i\bigr) 
  \II\{\, \lambda(t)g_i> R\,\}
  \lsim \log\varphi_0\bigl(\lambda(t)\bigr) |g|_\beta^\beta \, .
  \eqno{\equa{LDPCEqB}}
$$
Moreover, as was shown in the proof of Lemma \LDPBLemma, there exists 
a positive real number $c$ such that for any $n$ less than $\zeta V(t)$,
$$\eqalign{
  \sum_{0\leq i<n} \log\varphi_0\bigl( \lambda(t)g_i\bigr)
  \II\{\, \lambda(t) g_i\leq R\,\}
  &{}\leq c\sum_{0\leq i<n}\lambda(t)^2 g_i^2 \cr
  &{}\leq c\lambda(t)^2 \sum_{0\leq i<\zeta V(t)} g_i^2 \cr
  &{}=O\Bigl(\lambda(t)^2 {t^2\over V(t)}\Bigr)\, .\cr
  }
$$
As a function of $t$, this asymptotic upper bound is regularly varying
of index
$$
  {2\over \beta-1}+2-{1\over \gamma} = 2\alpha-{1\over\gamma} \, .
$$
The upper bound \LDPCEqB\ is regularly varying of index $\beta/(\beta-1)
=\alpha$. Since $\alpha\gamma$ is less than $1$, we see that $2\alpha-1/\gamma$
is less than $\alpha$, and, consequently, for $n$ less than $\zeta V(t)$,
$$\eqalign{
  \sum_{0\leq i<n}\log\varphi_0\bigl(\lambda(t)g_i\bigr)
  &{}\lsim \log\varphi_0\bigl(\lambda(t)\bigr) |g|_\beta^\beta \cr
  &{}\sim (\log\varphi_0)\circ m_0^\leftarrow (t)|g|_\beta^{-\alpha} \, .\cr
  }
$$
This implies that the exponent in the upper bound \LDPCEqA\ is 
asymptotically bounded by an equivalent of
$$
  -|g|_\beta^{-\alpha}m_0^\leftarrow(t)t+|g|_\beta^{-\alpha} \log\varphi_0\circ
  m_0^\leftarrow(t) 
  = -|g|_\beta^{-\alpha} I_0(t) \, .
$$
The result follows.\hfill\qed

\bigskip

We now prove a trivial lower bound.

\Lemma{\label{LDPD}
  For any positive $n$,
  $$
    \log\P\{\, S_n>t\,\} 
    \gsim \log\oF(t) \Bigl(\sum_{0\leq i<n}g_i^\beta\Bigr)^{-\alpha/\beta}
  $$
  as $t$ tends to infinity.
}

\bigskip

\Proof Let $x_i$ be $g_i^{1/(\alpha-1)}\bigm/\sum_{0\leq i<n}g_i^\beta$, so 
that $\sum_{0\leq i<n} g_ix_i=1$. We have
$$\eqalign{
  \log\P\{\, S_n>t\,\}
  &{}\geq \log\P\Bigl({\textstyle\bigcap}_{0\leq i<n} \{\, X_i>tx_i\,\}
             \Bigr) \cr
  \noalign{\vskip 3pt}
  &{}=\sum_{0\leq i<n}\log\oF(tx_i) \cr
  &{}\sim \log\oF(t) \sum_{0\leq i<n} x_i^\alpha \cr
  }
$$
as $t$ tends to infinity. The result follows upon calculating
$$
  \sum_{0\leq i<n} x_i^\alpha 
  = \Bigl(\sum_{0\leq i<n} g_i^\beta\Bigr)^{-\alpha/\beta} \, .
  \eqno{\qed}
$$

\bigskip

\noindent{\bf Proof of Theorem \ThCaseC}
{\it Lower bound.} Applying Proposition \LDPD, for any positive integer $n$,
$$\eqalign{
  \log\P\{\, M>t\,\}
  &{}\geq \log\P\{\, S_n>t\,\} \cr
  &{}\gsim \log\oF(t)\Bigl(\sum_{0\leq i<n}g_i^\beta\Bigr)^{-\alpha/\beta} \, .
   \cr
  }
$$
Consequently, as $t$ tends to infinity,
$$
  \log\P\{\, M>t\,\} \gsim \log\oF(t) |g|_\beta^{-\alpha} \, .
$$

\noindent{\it Upper bound.} We apply the remark following Proposition 
\MasterTh. In the present context, Proposition \LDPB\ shows that \HypUBLDP\
holds with $r_n=I_0(\g n)$ and $I(x)=|g|_\beta^{-\alpha}x^\alpha$. Note that
$\log\oF_0\sim\log\oF$ at infinity. Since  
Broniatowski and Fuchs' (1995) Theorem 3.1 implies that, under the assumption
of Theorem \ThCaseC, $-\log\oF_0\sim I_0$ at 
infinity, the function $r$ is regularly varying of index 
$\alpha\gamma$. Referring to Proposition \MasterTh, we see that 
$\theta=|g|_\beta^{-\alpha}$ for
$$
  \inf_{x\geq 0} x^{\alpha\gamma}(x^{-\gamma}+1)^\alpha
  = \inf_{x\geq 0}(1+x^\gamma)^\alpha
  = 1 \, .
$$
Since $\g{V(t)}\sim t$, 
$$
  -r\circ V(t)
  = -I_0(\g{V(t)})
  \sim -I_0(t)
  \sim \log\oF(t)
$$
as $t$ tends to infinity. Therefore, in view of this and Lemma \LDPC,
we see that condition \HypSmallN\ holds. This proves Theorem \ThCaseC\
when $\mu$ is $-1$.  The same scaling argument as in the end of the proof of
Theorem \ThCaseA\ allows for the extension to other values of
$\mu$.\hfill\qed

\bigskip



\section{Proof of Theorem \ThConditionalA.}
As for the proof of the results of section 2, we will first prove the
result when the mean $\mu$ is $-1$.
In the first two subsections, we prove assertion (i) and (ii) respectively.
Assertion (iii) requires a distinction between the cases of boundedness or 
divergence of $\max_{0\leq i\leq n} g_{i/n}$, and is proved, accordingly,
in the the third and fourth subsections. A scaling argument, developed
in the fifth subsection gives Theorem \ThConditionalA\ when the mean $\mu$
is arbitrary.

Throughout this section we will use the following obvious fact. Let $E_t$
be an event indexed by $t$. To prove that $\P(E_t\mid M>t)$ tends to $0$ as
$t$ tends to infinity, it suffices to prove that $\P(E_t)=o(\P\{\, M>t\,\})$
as $t$ tends to infinity; indeed, this follows from the definition of
conditional probability and monotonicity of measures.

\bigskip

\subsection{Proof of Theorem \ThConditionalA.i}
Assume that $\mu$ is $-1$. Assumptions of Proposition \MasterTh\ are
satisfied by virtue of Proposition \LDP\ and Lemma \ThASmallIndices.
From the second assertion of Proposition \MasterTh\ and equality \JJeq,
we deduce that $\calN_t$ converges to $\tau$ in probability as $t$
tends to infinity and conditionally on $M$ exceeding $t$. This is
assertion (i) when the mean is $-1$.

\bigskip

\subsection{Proof of Theorem \ThConditionalA.ii}
We assume that $\mu$ is $-1$.
Our next lemma is the analogue of Lemma \LDPLemma\ specialized to the context
of the proof of Theorem \ThConditionalA. Recall that $X_i$ has mean $-1$
for the time being, and that $\varphi_0$ is the moment generating function
of the centered random variable $X_i+1$.

\Lemma{\label{LemmaApproxE}
  Let $f$ be a continuous real-valued and bounded function 
  on $[\,0,1\,]\times \RR$. For any fixed $\lambda$,
  $$\displaylines{\qquad
      \limn n^{-1}\sum_{1\leq i\leq n} E\Bigl( f\Bigl({i\over n},X_i\Bigr)
      {e^{\lambda g_{n-i/n}(X_i+1)}\over \varphi_0(\lambda g_{n-i/n})}\Bigr)
    \hfill\cr\hfill
      {}=\int f(v,x) { e^{\lambda k_\gamma(v)(x+1)}
                       \over\varphi_0\bigl( \lambda k_\gamma(v)\bigr) }
      \II_{[0,1)}(v) \d (L\otimes F)(v,x) \, .
    \qquad\cr}
  $$
  }

\Proof Let $c$ be a number larger than $\limn \max_{0\leq i\leq n} g_{i/n}$.
Consider the function
$$
  \psi (v,y,x) 
  = f(v,x)
  {e^{\lambda (y\wedge c)(x+1)}\over\varphi_0 \bigl(\lambda(y\wedge c)\bigr)}
  \, .
$$
For $n$ large enough and with $\Gamma_n$ the measure defined prior to Lemma 
\GammaNConverges,
$$
  n^{-1}\sum_{1\leq i\leq n} \E \Bigl( f\Bigl({i\over n},X_i\Bigr)
  {e^{\lambda g_{n-i/n}(X_i+1)}\over\varphi_0\bigl(\lambda g_{n-i/n})}\Bigr)
  = \E\int \psi(v,y,X_1)\d{\mit\Gamma}_n(v,y) \, .
$$
For any fixed $x$ the function $\psi(v,y,x)$ is a continuous and bounded
function of $(v,y)$ in $[\,0,1\,]\times \RR$. By Lemma \GammaNConverges, the
sequence of functions
$$
  \psi_n(x)=\int \psi(v,y,x)\d{\mit\Gamma}_n(v,y)\, , \qquad n\geq 1\, ,
$$
converges pointwise to the function
$$
  \psi(x)=\int_0^1 \psi\bigl( u,k_\gamma(u),x\bigr) \d u \, .
$$
Since
$$
  |\psi (v,y,x)|
  \leq |f|_{[0,1]\times\RR} e^{\lambda c|x+1|}
  \Bigl|{1\over\varphi_0}\Bigr|_{[0,\lambda c]} \, ,
$$
the dominated convergence theorem implies that $\E\psi_n(X_1)$ tends 
to $\E\psi(X_1)$ as $n$ tends to infinity, which is what the lemma 
asserts.\hfill\qed

\bigskip

Recall that in section 2 we used the notation $\theta$ for 
$$
  \theta
  =-\limt V(t)^{-1}\log\P\{\, M>t\,\}
  \, .
$$
Considering the definition of $\tau$ in \HypTauUnique, that of $J_0$
in \JZeroDef, equality \JJeq, and how $\theta$ was obtained in \ThetaJZero,
$$\eqalignno{
  \theta
  &{}=\tau J_0(\tau^{-\gamma}+1) \cr
  &{}=\tau\sup_{\lambda>0}\Bigl( (\tau^{-\gamma}+1)\lambda
    -\int_0^1\log\varphi_0 \bigl(\lambda k_\gamma(u)\bigr)\d u \Bigr)\, .
  &\equa{EqTauJ} \cr
  }
$$
Since $m_0$ is onto the nonnegative half-line, the supremum in $\lambda$ in the
above formula is achieved for some value $A$. By considering the derivative
in $\lambda$, which must vanish at the maximizer $A$, we obtain
$$
  \tau^{-\gamma}+1
  = \int_0^1 k_\gamma(u) m_0\bigl(Ak_\gamma(u)\bigr)\d u \, .
  \eqno{\equa{ADef}}
$$
When $\mu$ is $-1$ as currently, we have $\varphi(\lambda)
=e^{-\lambda}\varphi_0(\lambda)$ and, consequently, $m=-1+m_0$. Therefore, the
definition of $A$ in \ADef\ matches that in \ADefA.

As will be apparent in the bound \fixedref{6.2.8} to come and in
its evaluation, the following result is strongly related to
Proposition \LDP\ if one takes $n$ to be about $\tau V(t)$ and $x$ to
be about $\tau$ in that proposition.

\Lemma{\label{LemmaApproxF}
  The following holds,
  $$\displaylines{
    \limeps\limsupt\sup_{n\,:\,{\textstyle |}{\ss n\over\ss V(t)}
                         -\tau{\textstyle |}<\epsilon}
    \Bigl| A \gamma {t+s_n\over g_nV(t)} -{1\over V(t)}\sum_{1\leq i\leq n}
    \log\varphi_0 (A g_{n-i/n})-\theta\Bigr| 
  \hfill\cr\hfill
  = 0 \, .\qquad}
  $$
}

\Proof Write $n=\nu V(t)$. Since $s$ is regularly varying, 
$s_n\sim \nu^\gamma t$. Moreover, \gEquiv\ shows that
$$
  g_n
  \sim \gamma \nu^{\gamma-1}{t\over V(t)} \, ,
$$
and those equivalences hold locally uniformly in $\nu$ thanks to the uniform
convergence Theorem (Bingham, Goldie and Teugels, 1989, Theorem 1.2.1).
In particular,
$$
  \gamma {t+s_n\over g_n V(t)}
  \sim {1+\nu^\gamma\over\nu^{\gamma-1}} \, ,
  \eqno{\equa{EqLemmaApproxFA}}
$$
as $t$ tends to infinity. Applying Lemma \LDPLemma, we also have
$$\eqalignno{
  {1\over V(t)}\sum_{1\leq i\leq n} \log\varphi_0 (A g_{n-i/n})
  &{}\sim {n\over V(t)} \int_0^1 \log\varphi_0\bigl(A k_\gamma(v)\bigr)\d v \cr
  &{}\sim \nu \int_0^1\log\varphi_0\bigl(A k_\gamma(v)\bigr)\d v \, ,
  &\equa{EqLemmaApproxFB}\cr}
$$
again locally uniformly in $\nu$ positive and as $t$ tends to infinity. 
Combining \EqLemmaApproxFA\ and \EqLemmaApproxFB, we obtain that
$$\displaylines{
    A \gamma{t+s_n\over g_n V(t)} 
    - {1\over V(t)}\sum_{1\leq i\leq n} \log\varphi_0(A g_{n-i/n})
  \hfill\cr\hfill
    {}=\nu \Bigl( A(\nu^{-\gamma}+1)
    -\int_0^1 \log\varphi_0\bigl(A k_\gamma(v)\bigr)\d v \Bigr) + o(1)
  \cr\hfill
  \equa{EqLemmaApproxFC}
  \cr}
$$
as $t$ tends to infinity. When $\nu$ is $\tau$, equality \EqTauJ\ shows that 
the right hand side in \EqLemmaApproxFC\ is $\theta$.
The result follows from the continuity in $\nu$ of the function involved in
the right hand side of \EqLemmaApproxFC.\hfill\qed

\bigskip

We can now prove the second assertion of Theorem \ThConditionalA. 
Let $f$ be a continuous function
supported by a vertical strip of the right half-space. 
Whenever $\nu$ is a 
measure on the right half-space, we write $\nu f$ for $\int f\d\nu$.
Let $\epsilon$ be a positive real number. Assume that we have proved that 
for any real number $h$ greater than $\calM f$,
$$
  \limt \P\{\, \calM_tf>h\mid M>t\,\} =0 \, . 
  \eqno{\equa{iiiA}}
$$
If $h$ is less then $\calM f$, then applying the above relation to $-f$ and 
$-h$, we see that the conditional probability of $\calM_t f<h$ given $M$
exceeds $t$ tends to $0$ as $t$ tends to infinity. We then conclude that
$$
  \limt\P\{\, |(\calM_t-\calM)f|>\epsilon\mid M>t\,\} = 0 \, .
$$
Thus, as $t$ tends to infinity, $\calM_tf$ converges in probability to
$\calM f$ conditionally on $M$ exceeding $t$. Since $f$ is arbitrary,
this shows that $\calM_t$ converges to $\calM$ in probability, under
the conditional probability that $M$ exceeds $t$. This would prove the
second assertion of Theorem \ThConditionalA, and therefore, it
suffices to prove \iiiA, which we do now.

The proof of the second assertion of Proposition \MasterTh\ shows 
that for any positive $\epsilon$,
$$\displaylines{\quad
  \P\{\, \calM_t f>h\, ;\, M>t\,\}
  \hfill\cr\hfill
  \leq \sum_{n\,:\, {\textstyle|}{\ss n\over\ss V(t)}-\tau{\textstyle|}\leq
    \epsilon} \P\{\, \calM_tf>h\,;\, S_n^0>t+s_n\,\} + o(\P\{\, M>t\,\})
  \quad\equa{iiiB}\cr}
$$
as $t$ tends to infinity. 

The basic inequality for our proof is the exponential form of Markov's,
which implies that for any positive $\lambda$,
$$\eqalignno{
  {}&\ \hskip -20pt\log\P\{\, \calM_t f>h \,;\, S_n^0>t+s_n\,\}\cr
  {}\leq{}&\log\P\Bigl\{\, \lambda V(t)\calM_t f+A\gamma {S_n^0\over g_n}
    >\lambda V(t) h+A \gamma{t+s_n\over g_n}\,\Bigr\} \cr
  {}\leq{}&-V(t)\Bigl( \lambda h+A\gamma {t+s_n\over g_nV(t)}
    -{1\over V(t)}\log\E\exp\Bigl( \lambda V(t)\calM_t f +A\gamma 
    {S_n^0\over g_n}\Bigr)\Bigr) \, .\cr
  & &\equa{iiiE}
  }
$$
The remainder of the proof is somewhat technical, but the next few sentences
show that it is very simple in essence. If one looks at the classical
Chernoff estimate, one sees that the minimizer in the rate function tends to
$0$ as one considers deviations nearing the mean. In our case, \iiiE\ is to be
considered when $h$ is close to the anticipated limit $\calM f$. Therefore, we
anticipate that we may take $\lambda$ very small. If this is so, the bound
can be linearized in $\lambda$. The linear term in $\lambda$ will be positive,
while the term measuring the deviation of $S_n^0>t+s_n$ should give a 
contribution very close to $\theta$. So, the addition of the linear term
in $\lambda$ to the term near $\theta$ should give a term greater 
than $\theta$, which is all that we need.

To proceed rigorously,
we define some small ---~arguably, bewildering~--- constants. Let $\delta$
be a positive real number less than $1$ such that $h>(1+2\delta)\calM f$. 
Let $\lambda$ be positive and small enough so that
$$
  \lambda |f|_{[0,\infty)\times\RR}
  <\sup\{\, x\,:\, e^x<1+(1+\delta)x\,\}\wedge (1+\delta)^{-1} \, .
  \eqno{\equa{iiiC}}
$$
Next, let $\eta$ be small enough so that $\lambda \bigl( h-(1+2\delta)\calM f
\bigr)>3\eta$. Finally, using Lemma \LemmaApproxF, let $\epsilon$ be a 
positive real number so that
$$\displaylines{
    \limsupt\sup_{n\,:\,{\textstyle |}{\ss n\over\ss V(t)}-\tau{\textstyle |}
    <\epsilon}
    \Bigl| A \gamma{t+s_n\over g_nV(t)} -{1\over V(t)}\sum_{1\leq i\leq n}
    \log\varphi_0 (A g_{n-i/n})-\theta\Bigr|  
  \hfill\cr\hfill
  {}< \eta \, . \qquad \equa{iiiD}\cr}
$$

To evaluate the upper bound \iiiE, we first bound the term containing an 
expectation. Given the constraint \iiiC\ on $\lambda$,
$$
  e^{\lambda f(i/V(t),X_i)} 
  \leq  1+(1+\delta)\lambda f\Bigl({i\over V(t)},X_i\Bigr) \, .
$$
Recall that $X_i$ is of mean $-1$ currently.
Since
$$\displaylines{
  \lambda V(t)\calM_t f+A \gamma S^0_n/g_n
  \hfill\cr\hfill
  {}=
  \sum_{1\leq i\leq n} \Bigl(\lambda f\Bigl({i\over V(t)},X_i\Bigr)
                             +A g_{n-i/n}(X_i+1)\Bigr)
  +\sum_{i>n}\lambda f\Bigl({i\over V(t)},X_i\Bigr) \, ,
  \cr}
$$
the term
$\E\exp\Bigl(\lambda V(t)\calM_tf +A \gamma{\ds S_n^0\over\ds g_n}\Bigr)$
in \iiiE\ is at most
$$\displaylines{\qquad
  \prod_{1\leq i\leq n}\E\Bigl( 1+(1+\delta)\lambda 
           f\Bigl({i\over V(t)}, X_i\Bigr)\Bigr)
           e^{A g_{n-i/n}(X_i+1)}
  \hfill\cr\hfill
  \prod_{i>n} \E\Bigl(1+(1+\delta)\lambda f\Bigl({i\over V(t)},X_i\Bigr)\Bigr)
  \, .
  \qquad\equa{iiiF}\cr}
$$
Note the inequality $\log(a+b)\leq \log a+b/a$, valid for any positive $a$ 
and any $b$ larger than $-a$. To apply this inequality with
$$
  a=\E e^{Ag_{n-i/n}(X_i+1)}
$$
and
$$
  b=(1+\delta) \lambda \E f\Bigl({i\over V(t)},X_i\Bigr) e^{Ag_{n-i/n}(X_i+1)}
  \, ,
$$
we first observe that
$$
  |b|
  \leq (1+\delta) \lambda |f|_{[0,\infty)\times\RR} \E e^{Ag_{n-i/n}(X_i+1)}
$$ 
and \iiiC\ ensures that $|b|$ is less than $a$. Therefore, $a+b$ is positive. 
We then have, referring to the first product of \iiiF,
$$\eqalign{
  &\hskip-12pt
   \log\E\Bigl( 1+(1+\delta)\lambda f\Bigl({i\over V(t)},X_i\Bigr) \Bigr)
   e^{A g_{n-i/n}(X_i+1)} \cr
  &{}\leq \log\E e^{A g_{n-i/n}(X_i+1)} \cr
  &\qquad{}+(1+\delta) \lambda 
   { \E f\bigl(i/V(t),X_i)e^{A g_{n-i/n}(X_i+1)}
     \over
     \E e^{A g_{n-i/n}(X_i+1)} } \cr
  &{}=\log\varphi_0 (A g_{n-i/n}) +(1+\delta)\lambda
    \E \Bigl( f\Bigl({i\over V(t)},X_i\Bigr) { e^{A g_{n-i/n}(X_i+1)}
                              \over
                              \varphi_0(A g_{n-i/n}) }
       \Bigr) \, .\cr}
$$
Consequently, using the inequality $\log (1+x)\leq x$ to handle the second
product in the upper bound \iiiF, we see that the logarithm of \iiiF\ is at
most
$$\displaylines{
  \sum_{1\leq i\leq n}\hskip-2pt \log\varphi_0 (A g_{n-i/n})
  +(1+\delta)\lambda\hskip-3pt\sum_{1\leq i\leq n} \hskip-2pt
  \E\Bigl( f\Bigl({i\over V(t)},X_i\Bigr) 
   { e^{A g_{n-i/n}(X_i+1)}
    \over
    \varphi_0(A g_{n-i/n}) }\Bigr)
  \hfill\cr\hfill
  {}+(1+\delta)\lambda\sum_{i>n} \E f\Bigl( {i\over V(t)},X_i\Bigr) \, .
  \qquad\cr}
$$
Referring to the upper bound \iiiE,
$$
  \lambda h +A\gamma {t+s_n\over g_n V(t)} 
  -{1\over V(t)}
  \log\E\exp\Bigl( \lambda V(t)\calM_t f+A \gamma{S^0_n\over g_n}\Bigr)
$$
is then at least
$$\eqalignno{
  A\gamma {t+s_n\over g_n V(t)} 
  &{}-{1\over V(t)}\sum_{1\leq i\leq n}
  \log\varphi_0(A g_{n-i/n}) +\lambda h \cr
  &{}-(1+\delta) \lambda {n\over V(t)} {1\over n} \sum_{1\leq i\leq n} 
     \E\Bigl( f\Bigl( {i\over V(t)},X_i\Bigr) 
     { e^{A g_{n-i/n}(X_i+1)}
     \over
     \varphi_0(A g_{n-i/n}) } \Bigr)\cr
  &{}-(1+\delta)\lambda {n\over V(t)} {1\over n}
   \sum_{i>n}\E f\Bigl({i\over V(t)},X_i\Bigr) \, .
  &\equa{iiiG}\cr
  }
$$
Define $\nu$ as $n/V(t)$. Using \iiiD, Lemma \LemmaApproxE, the equality
$\varphi_0(\lambda)=e^\lambda\varphi(\lambda)$ valid here since $\mu$ is $-1$,
we obtain that \iiiE\ is at most the exponential of $-V(t)$ times
$$\eqalignno{
    \theta&{}-\eta +\lambda h\cr
    &{}-(1+\delta)\lambda\nu\int f(v \nu,x)
      { e^{A \gamma(1-v)^{\gamma-1}(x+1)} 
        \over 
        \varphi_0\bigl(A\gamma(1-v)^{\gamma-1}\bigr) 
      }
      \II_{[0,1)}(v) \d(L\otimes F)(v,x) \cr
    &{}-(1+\delta) \lambda\nu \int f(v\nu,x)\II_{[1,\infty)} (v)
      \d(L\otimes F)(v,x) \cr
    &\qquad{}= \theta-\eta +\lambda \Bigl( h-(1+\delta){\nu\over\tau}\int 
      f\Bigl({\nu\over\tau}v,x\Bigr) \d\calM (v,x)\Bigr) \, .
    &\equa{iiiH}\cr}
$$
If $\epsilon$ is small enough so that $\nu/\tau$ is close enough to $1$, then
$$
  \Bigl|\int {\nu\over\tau}f\Bigl({\nu\over\tau}v,x\Bigr) \d\calM (v,x)
        -\calM f\Bigr|
  <\eta/\lambda
$$
and \iiiH\ is at least
$$
  \theta-\eta +\lambda \bigl( h-(1+\delta)\calM f) -\eta \, ,
$$
which, by our choice of $\eta$ is greater than $\theta+(1-\delta)\eta$. Hence
$$
  \log\P\{\, \calM_t f>h\,;\, S^0_n>t+s_n\,\} 
  \lsim -V(t)\bigl(\theta+(1-\delta)\eta\bigr)
$$
as $t$ tends to infinity. Since $V$ is regularly varying, \iiiB\ shows that
\iiiA\ holds, and this proves assertion (ii) of Theorem \ThConditionalA\
when $\mu$ is $-1$.

\bigskip

\subsection{Proof of Theorem \ThConditionalA.iii when 
\poorBold{$\max_{0\leq i\leq n} g_{i/n}$} is ultimately bounded}
In essence, the proof consists in writing the process $\calS_t$ as a functional
of $\calM_t$ and showing that the convergence of $\calM_t$ to $\calM$ implies
that of the functional of $\calM_t$ to the functional of $\calM$. The main
difficulty is that the functional is not continuous with respect to our
topology on measures. This forces us to develop various approximation results
to show that $\calS_t$ is approximable by a well behaved functional 
of $\calM_t$.

To proceed, for any measure $\nu$ on the right half-space for which the 
integrals
$$
  \int \II_{[0,\lambda)}(v) (\lambda -v)^{\gamma-1} |x| \d\nu(v,x) \, ,
  \qquad \lambda >0 \, ,
$$
are finite, we define the functional $\fS$ of $\nu$ evaluated at $\lambda$
by
$$
  \fS (\nu)(\lambda)= 
  \int\II_{[0,\lambda)}(v) \gamma(\lambda -v)^{\gamma-1} x \d\nu(v,x) \, ,
$$
with the convention that $\fS(\nu)(0)$ is $0$.

For any function $f$ defined on some interval $[\,a,b\,]$ we write
$$
  |f|_{[a,b]}=\sup_{a\leq x\leq b} |f(x)|
$$
for its supremum norm over that interval.

Our first lemma shows that given that $M$ exceeeds a large threshold $t$, 
the process $\St$ is well approximated by $\fS (\calM_t)$ locally
uniformly.

\Lemma{\label{LemmaApproxA}
  For any positive $\Lambda$ and $\epsilon$,
  $$
    \limt \P\bigl\{\, |\St-\fS(\calM_t)|_{[0,\Lambda]}>\epsilon \bigm| M>t\,
            \bigr\} = 0 \, .
  $$
}

\Proof Consider the difference $\Delta_t=|\St-\fS(\calM_t)|_{[0,\Lambda]}$.
To analyse it, we rewrite $\St(\lambda)$ as
$$\eqalignno{
  \St (\lambda)
  &{}={1\over t} \sum_{1\leq i\leq \lambda V(t)} 
    g_{\lfloor \lambda V(t)\rfloor-i} X_i \cr
  &{}=\int \II_{[0,\lambda)}(v) {V(t)\over t} g_{\lfloor\lambda V(t)\rfloor-
     \lfloor vV(t)\rfloor} x\d\calM_t(v,x) \cr
  &\qquad\qquad\qquad\qquad\qquad 
   {}+{1\over t} \II_\NN \bigl(\lambda V(t)\bigr) g_0 
     X_{\lfloor\lambda V(t)\rfloor} \, . 
  &\equa{EqLemmaApproxAA} \cr
  }
$$
From this expression and the following consequence of \gEquiv,
$$
  {V(t)\over t} g_{\lfloor\lambda V(t)\rfloor-\lfloor vV(t)\rfloor}
  \sim \gamma (\lambda-v)^{\gamma-1} \, ,
  \eqno{\equa{EqLemmaApproxAAA}}
$$
the result appears natural, though not proved yet. The proof has four steps.

We fix a positive $\Lambda$ and
we consider a positive real number $\eta$.

\noindent {\it Step 1.} Let $\delta$ be a positive real number
and define $\Delta_{t,1}(\lambda)$ as
$$
  \int \II_{[0,(\lambda-\delta)_+)}(v) 
  \Bigl| {V(t)\over t}  g_{\lfloor\lambda V(t)\rfloor-\lfloor vV(t)\rfloor}
  -\gamma(\lambda-v)^{\gamma-1}\Bigr| |x| \d\calM_t (v,x) \, .
$$
The asymptotic equivalence \EqLemmaApproxAAA\ holds uniformly in the 
range of $\lambda$ and $v$ such 
that $0\leq v<v+\delta<\lambda<\Lambda$.
Consequently, for $t$ large enough and uniformly in $\lambda$ 
in $[\,\delta,\Lambda\,]$, 
$$
  \Delta_{t,1}(\lambda)
  \leq \eta\int\II\{\, 0\leq v<v+\delta<\lambda\,\} 
  \gamma (\lambda -v)^{\gamma-1} |x|\d \calM_t(v,x)\, .
$$
Since $(\lambda-v)^{\gamma-1}\leq\Lambda^{\gamma-1}$ in that range, we further
obtain the upper bound, independent of $\lambda$,
$$\eqalign{
  \Delta_{t,1}(\lambda)
  &{}\leq \eta\gamma \Lambda^{\gamma-1} 
    \int\II\{\, 0<u<u+\delta <\Lambda\,\} |x|\d\calM_t(u,x) \cr
  &{}\leq \eta \gamma\Lambda^\gamma {1\over \Lambda V(t)} 
   \sum_{1\leq u\leq \Lambda V(t)} |X_i| \, .\cr
  }
$$

We take $\eta$ small enough so that 
$\epsilon/\eta\gamma\Lambda^\gamma$ exceeds the mean of $|X_1|$.
Since the moment generating function of $|X_1|$ is finite in a neighborhood
of the origin, we introduce the Cram\'er function associated to the 
distribution of $|X_1|$, 
$$
  I_*(x)=\sup_{s>0} sx-\log \E e^{s|X_1|} \, .
$$
The Chernoff bound implies,
$$
  \P\{\, |\Delta_{t,1}|_{[\delta,\Lambda]}>\epsilon\,\}
  \leq \exp\Bigl( -\lfloor\Lambda V(t)\rfloor 
  I_*\Bigl( {\epsilon\over\eta\gamma\Lambda^\gamma}\Bigr)\Bigr)
  \, .
  \eqno{\equa{EqLemmaApproxAB}}
$$
Since $ I_*(\epsilon/\gamma\Lambda^\gamma\eta)$ tends to infinity 
as $\eta$ tends to $0$, we
can choose $\eta$ small enough to guarantee that the upper bound 
\EqLemmaApproxAB\ is negligible compared to the probability that $M$ exceeds 
$t$ by virtue of Theorem \ThCaseA. Consequently,
$$
  \limt \P\bigl\{\, |\Delta_{t,1}|_{[\delta,\Lambda]}>\epsilon \bigm| M>t\,\}
  = 0 \, .
$$

\noindent{\it Step 2.} We now consider
$$
  \Delta_{t,2}(\lambda)=
  \Bigl| \int\II_{[(\lambda-\delta)_+,\lambda)}(v) {V(t)\over t} 
  g_{\lfloor\lambda V(t)\rfloor-\lfloor vV(t)\rfloor} x \d\calM_t(v,x) \Bigr|\, .
$$
Let $c$ be a positive number such that $\max_{0\leq i\leq n}g_i \leq cg_n$
for any $n$ large enough. For $t$ large enough and uniformly 
in $\lambda$ in $[\,0,\Lambda\,]$,
$$\eqalign{
  \Delta_{t,2}(\lambda)
  &{}\leq c {V(t)\over t} g_{\lfloor\delta V(t)\rfloor}
    \int\II_{[(\lambda-\delta)_+,\lambda)}(v) |x|\d\calM_t(v,x) \cr
  &{}\leq 2c\gamma\delta^{\gamma-1} {1\over V(t)} \sum_{(\lambda-\delta)_+V(t)
    \leq i<\lambda V(t)} |X_i| \, . \cr
  }
$$
Therefore,
$$
  |\Delta_{t,2}|_{[0,\Lambda]} 
  \leq 2c\gamma\delta^\gamma\max_{1\leq j<\Lambda V(t)}
  {1\over\delta V(t)}\sum_{j\leq i<j+\delta V(t)} |X_i| \, .
  \eqno{\equa{EqLemmaApproxABB}}
$$
Using Bonferroni's inequality and then Chernoff's, this implies that for 
$t$ large enough,
$$\eqalignno{
  \P\{\, |\Delta_{t,2}|_{[0,\Lambda]} >\epsilon\,\}
  &{}\leq \Lambda V(t) 
    \P\Bigl\{\, \sum_{1\leq i\leq\delta V(t)} |X_i| >
                { \delta V(t)\over 2c\gamma\delta^\gamma}\epsilon \,\Bigr\} 
    \cr
  &{}\leq\Lambda V(t)
    \exp\Bigl( -\lfloor\delta V(t)\rfloor
               I_*\Bigl( {\epsilon\over 2c\gamma\delta^\gamma}\Bigr)\Bigr) \, .
  &\equa{EqLemmaApproxAC}\cr
  }
$$
Since $I_*\gg \Id$ at infinity, taking $\delta$ small enough ensures that
the upper bound \EqLemmaApproxAC\ is negligible compared to the probability
that $M$ exceeds $t$.

\noindent{\it Step 3.} We now consider
$$
  \Delta_{t,3}(\lambda)
  =\Bigl|\int\II_{[(\lambda-\delta)_+,\lambda)}(v)
   \gamma(\lambda-v)^{\gamma-1} x\d\calM_t(v,x)\Bigr| \, .
$$
This is at most
$$
  \gamma\delta^{\gamma-1} {1\over V(t)} 
  \sum_{(\lambda-\delta)_+V(t)\leq i<\lambda V(t)} |X_i| \, .
$$
Comparing with \EqLemmaApproxABB, we deduce from the previous step that for 
any $\delta$ small enough,
$$
  \limt\P\bigl\{\, |\Delta_{t,3}|_{[0,\Lambda]}>\epsilon \bigm| M>t\,\bigr\}
  = 0 \, .
$$

\noindent{\it Step 4.} Let
$$
  \Delta_{t,4}(\lambda)
  =t^{-1} \II_\NN\bigl(\lambda V(t)\bigr) g_0 
  X_{\lfloor \lambda V(t)\rfloor} \, .
$$
We see that
$$
  |\Delta_{t,4}|_{[0,\Lambda]}
  =t^{-1} g_0\max_{1\leq i\leq\Lambda V(t)} X_i \, .
$$
Thus, Bonferroni's inequality yields
$$
  \P\{\, |\Delta_{t,4}|_{[0,\Lambda]} >\epsilon\,\}
  \leq \Lambda V(t) \oF(t\epsilon/g_0) \, .
$$
To prove that this upper bound is negligible compared to the probability
that $M$ exceeds $t$, it suffices to show that for any $\epsilon$ and
$\eta$ positive, $\limt \log\oF(\epsilon t)+\eta V(t)=-\infty$. This limit
holds for the following reason. Firstly, the finiteness of the moment
generating function on the nonnegative half-line implies that 
$\limt \log\oF(\epsilon t)/t=-\infty$ for any positive $\epsilon$. Secondly, 
Lemma \UGrowth\ implies that
$\limsup_{t\to\infty} V(t)/t$ is finite. Thus, $\limt t^{-1}\bigl(
\log\oF(\epsilon t)+\eta V(t)\bigr)=-\infty$ for any positive $\epsilon$ and
$\eta$, which is more than what we needed.

\noindent{\it Conclusion.} Since $\Delta_t$ is at most the sum
$\Delta_{t,1}+\Delta_{t,2}+\Delta_{t,3}+\Delta_{t,4}$, the result follows
from the four steps, Bonferroni's inequality and the fact that $\epsilon$
is arbitrary.\hfill\qed

\bigskip

The functional $\fS$ is not well behaved with respect to weak$*$ or vague
convergence, because it integrates a function which is both unbounded in $x$
and discontinuous in $v$. 
Those are classical problems which arise in large deviations theory when
one wants to use a so-called contraction principle, and the remedy is often
to use a truncation and a smoothing ---~in a context closely related to ours, 
see Bahadur (1971), Groeneboom, Oosterhoff and Ruymgaart (1979), and 
Hoadley (1967) for the truncation argument.
Other approaches, such as Ganesh and O'Connell's (2002), could likely be used
as well. To setup the truncation argument, let $b$ be a positive real number
and define, now for any measure $\nu$,
$$
  \fS(\nu,b)(\lambda)
  =\int\II_{[0,\lambda)}(v) \gamma (\lambda-v)^{\gamma-1} 
    \sign(x)(|x|\wedge b) \d \nu(v,x) \, .
$$
Note that for any positive $\Lambda$,
$$
  |\fS(\nu,b)-\fS(\nu)|_{[0,\Lambda]}
  \leq\gamma \Lambda^{\gamma-1}\int \II_{[0,\Lambda)}(v)(|x|-b)_+ \d\nu(v,x) 
  \, .
  \eqno{\equa{EqBoundTruncation}}
$$
Our next lemma shows that provided $b$ is large enough, $\fS(\calM_t,b)$ is
close to $\fS(\calM_t)$ in conditional probability given that $M$ exceeds a
large $t$.

\Lemma{\label{LemmaApproxB}
  For any $\Lambda$ and $\epsilon$ positive,
  $$
    \lim_{b\to\infty} \limsupt 
    P\bigl\{\, |\fS(\calM_t)-\fS(\calM_t,b)|_{[0,\Lambda]}>\epsilon \bigm| 
    M>t\,\bigr\} = 0 \, .
  $$
}

\Proof The upper bound in \EqBoundTruncation\ with $\calM_t$ substituted for
$\nu$ is 
$$
  \gamma\Lambda^{\gamma-1}{1\over V(t)} 
    \sum_{1\leq i\leq\Lambda V(t)} (|X_i|-b)_+ \, .
$$
Therefore, using the exponential Markov inequality, for any positive $a$,
$$\displaylines{\quad
  \P\bigl\{\, |\fS(\calM_t)-\fS(\calM_t,b)|_{[0,\Lambda]} >\epsilon \,\bigr\}
  \hfill\cr\hfill
  \leq \exp\Bigl( -\lfloor\Lambda V(t)\rfloor
                  \Bigl( a {\epsilon\over\gamma\Lambda^\gamma}-
                         \log \E e^{a(|X_i|-b)_+}\Bigr) 
           \Bigr) \, .
  \quad\equa{EqLemmaApproxBA}\cr
  }
$$
By dominated convergence with dominating function $e^{a|X_i|}$, for any fixed
$a$,
$$
  \lim_{b\to\infty} \log \E e^{a(|X_i|-b)_+} = 0 \, .
$$
Therefore, taking $a$ such that $a\epsilon/\gamma\Lambda^\gamma$ is
large enough, we obtain that, provided $b$ is large enough, the right
hand side of \EqLemmaApproxBA\ is negligible, as $t$ tends to
infinity, compared to the probability that $M$ exceeeds $t$.  The result
follows.\hfill\qed

\bigskip

Recall that the measure $\calM$ was defined in \MDef.  Our next lemma
shows that the deterministic functions $\fS(\calM,b)$ and $\fS(\calM)$
are close provided $b$ is large enough.

\Lemma{\label{LemmaApproxC}
  For any positive $\Lambda$,
  $$
    \lim_{b\to\infty} |\fS(\calM,b)-\fS(\calM)|_{[0,\Lambda]}=0\, .
  $$
}

\Proof
Since $(|x|-b)_+$ is at most $|x|$ and the function $\II_{[0,\Lambda]}(v)|x|$
is $\calM$-integrable, the dominated convergence theorem ensures that,
after substituting $\calM$ for $\nu$ in \EqBoundTruncation, the right hand
side of \EqBoundTruncation\ tends to $0$ as $b$ tends to infinity.\hfill\qed

\bigskip

We now calculate $\fS(\calM)$, showing that it is equal to the function $\calS$
defined in \SDef\ and involved in Theorem \ThConditionalA.

\Lemma{\label{LemmaApproxD}
  For any positive $\lambda$,
  $$
    \fS(\calM)(\lambda)
    =\int_0^\lambda
    \gamma (\lambda-v)^{\gamma-1} m\bigl( Ak_\lambda(v/\tau)\bigr)
    \d v \, .
  $$
}

\Proof It follows from the identity
$$
  \int x {\exp \bigl( Ak_\gamma(v/\tau)x \bigr)
         \over
           \varphi\bigl(A k_\gamma(v/\tau)\bigr)
         }\d F(x)
   = m\bigl( A k_\gamma(v/\tau)\bigr)\, .
  \eqno{\qed}
$$

\bigskip

We now consider the modulus of continuity of $\St$ at $\lambda$,
$$
  \omega_{t,\delta}(\lambda)
  = \sup\{\, |\St(\lambda+v)-\St(\lambda)|\,:\, |v|\leq\delta\,\} \, ,
  \qquad \delta>0 \, .
$$
Our next lemma shows that $\St$ is very likely to be nearly uniformly 
continuous  when $M$ exceeds a large threshold $t$.

\Lemma{\label{LemmaApproxDB}
  For any positive $\epsilon$ and $\Lambda$,
  $$
    \lim_{\delta\to 0} \limsupt
    \P\{\, |\omega_{t,\delta}|_{[0,\Lambda]}>\epsilon\mid M>t\,\} = 0 \, .
  $$
}

\Proof For any positive $\delta$, define
$$
  \omega_{t,\delta,b}(\lambda)
  =\sup\{\, |\fS(\calM_t,b)(\lambda+v)-\fS(\calM_t,b)(\lambda)|\,:\, 
            |v|\leq\delta\,\} \, .
$$
Since
$$
  |\omega_{t,\delta}|_{[0,\Lambda]}
  \leq |\omega_{t,\delta,b}|_{[0,\Lambda]} 
       + 2|\St-\fS(\calM_t,b)|_{[0,\Lambda]} \, ,
$$
Lemmas \LemmaApproxA\ and \LemmaApproxB\ show that it suffices to prove that
for any $b$,
$$
  \lim_{\delta\to 0}\limsupt\P\{\, |\omega_{t,\delta,b}|_{[0,\Lambda]}>\epsilon
  \mid M>t\,\} = 0 \, .
$$
Let $\lambda_1$ and $\lambda_2$ be two positive real numbers, with $\lambda_1
<\lambda_2\leq\lambda_1+\delta$ and $\lambda_2\leq\Lambda$. 
We bound $|\fS(\calM_t,b)(\lambda_2)-\fS(\calM_t,b)(\lambda_1)|$ as the sum of
$$
  \gamma b \int\bigl( |\II_{[0,\lambda_2)}(v)-\II_{[0,\lambda_1)}|(v)\bigr) 
  (\lambda_2-v)^{\gamma-1}\d\calM_t(v,x)
  \eqno{\equa{EqLemmaApproxDBA}}
$$
and
$$
  \gamma b \int\II_{[0,\lambda_1)}(v) \bigl( (\lambda_2-v)^{\gamma-1}
  - (\lambda_1-v)^{\gamma-1} \bigr) \d\calM_t(v,x) \, .
  \eqno{\equa{EqLemmaApproxDBB}}
$$
For any $t$ large enough, \EqLemmaApproxDBA\ is at most
$$\eqalign{
  \gamma b \Lambda^{\gamma-1}
  \calM_t\bigl([\lambda_1,\lambda_2)\times\RR)\bigr)
  &{}\leq \gamma b \Lambda^{\gamma-1} 
   \Bigl( \lambda_2-\lambda_1+{1\over V(t)}\Bigr) \cr
  &{}\leq \gamma b \Lambda^{\gamma-1} 2\delta \, . \cr
  }
$$
The second term of \EqLemmaApproxDBB\ is at most the following function
evaluated at $t$,
$$
  \gamma b {1\over V^\gamma} \sum_{1\leq i\leq \lambda_1 V}
  \Bigl( \bigl( \lambda_2 V-i\bigr)^{\gamma-1} 
    -\bigl(\lambda_1 V-i\bigr)^{\gamma-1}\Bigr) \, .
  \eqno{\equa{EqLemmaApproxDBC}}
$$
Using the comparison of a sum with an integral, namely that for any positive
real numbers $0<a<b$,
$$
  \gamma\sum_{1\leq i\leq a} (b-i)^{\gamma-1}
  \cases{\leq \gamma\int_0^a (b-u)^{\gamma-1}\d u 
          \leq b^\gamma\, , \cr
         \geq \gamma\int_1^{\lfloor a\rfloor} (b-u)^{\gamma-1}\d u
           = (b-1)_+^\gamma - (b-\lfloor a\rfloor)^\gamma \, ,\cr
        }
$$
we see that $\gamma$ times the sum involved in \EqLemmaApproxDBC\ is at most
$$\displaylines{\qquad
  \bigl(\lambda_2V\bigr)^\gamma -\bigl(\lambda_1 V-1\bigr)_+^\gamma
  +\bigl(\lambda_1 V-\lfloor\lambda_1 V\rfloor\bigr)^\gamma
  \hfill\cr\hfill
  {}\leq V^\gamma
  \Bigl( \lambda_2^\gamma -\Bigl(\lambda_1-{1\over V}\Bigr)_+^\gamma
  +{1\over V^\gamma }\Bigr) \, .
  \qquad\cr}
$$
Since the function $\Id^\gamma$ is locally uniformly continuous on the 
nonnegative half-line, this implies that \EqLemmaApproxDBC\ can be made
arbitrarily small by taking $t$ large enough and $\delta$ small enough. 
This proves the lemma.\hfill\qed

\bigskip

Next, we setup the smoothing procedure which will allow us to approximate
the functional $\fS(\cdot,b)$ by a well behaved one. Define the function
$$
  \Ieps(v) 
  = \cases{ 1 & if $0\leq v\leq 1-\epsilon$, \cr
            (1-v)/\epsilon & if $1-\epsilon \leq v\leq 1$,\cr
            0 & if $v\geq 1$. \cr}
$$
This function is continuous, coincides with $\II_{[0,1)}$ on the complement
of $(1-\epsilon,1)$ and moreover, $0\leq \II_{[0,1)}-\Ieps\leq 1$.
Define the functional
$$
  \fS_\epsilon (\nu,b)(\lambda)
  =\int\Ieps \Bigl({v\over\lambda}\Bigr)
  \gamma(\lambda-v)^{\gamma-1} \sign (x) (|x|\wedge b) \d\nu(v,x) \, .
$$
The following shows that $\fS(\cdot,b)$ is well approximated by 
$\fS_\epsilon(\cdot,b)$ for the measures of interest to us.

\Lemma{\label{LemmaApproxDA}
  For any positive $\Lambda$ and any $t$ large enough, both 
  $|\fS_\epsilon(\calM_t,b)-\fS(\calM_t,b)|_{[0,\Lambda]}$ and
  $|\fS_\epsilon(\calM,b)-\fS(\calM,b)|_{[0,\Lambda]}$ are bounded by
  $2\gamma\Lambda^\gamma b\epsilon$.
}

\bigskip

\Proof Let $\nu$ be a $\sigma$-finite measure on the right half-space and
let $\lambda$ be nonnegative and at most $\Lambda$. Since 
$0\leq \II_{[0,1)}-\Ieps\leq\II_{[1-\epsilon,1)}$,
$$
  \bigl|\bigl(\fS (\nu,b)-\fS_\epsilon(\nu,b)\bigr)\bigr| (\lambda)
  \leq \gamma \Lambda^{\gamma-1}b 
  \nu\{\, (v,x)\,:\, 1-\epsilon\leq v/\lambda\leq 1\,\} \, .
$$
If $\nu$ is $\calM$, its first marginal measure is the Lebesgue measure; 
then $\nu\{\, (v,x)\,:\, 1-\epsilon\leq v/\lambda\leq 1\,\}$
is equal to $\lambda\epsilon$, while, if $\nu$ is $\calM_t$, it is equal to
$$ 
  {1\over V(t)} \sharp\{\, i\,:\, (1-\epsilon)\lambda V(t) 
   \leq i\leq\lambda V(t)\,\}
  \leq {1\over V(t)} \bigl(\lambda\epsilon V(t)+1\bigr) \, .
$$
The result follows.\hfill\qed

\bigskip

We now prove assertion (iii) of Theorem \ThConditionalA.  Combining
Lemmas \LemmaApproxA, \LemmaApproxB, \LemmaApproxDA\ and
\LemmaApproxC, we deduce that as $t$ tends to infinity and for any
$\lambda$ fixed, $\St(\lambda)$ converges in
probability to $\fS(\calM)(\lambda)$ conditionally on $M$ exceeds $t$.
Lemma \LemmaApproxDB\ turns this pointwise convergence to a locally
uniform one. The result follows from Lemma \LemmaApproxD.

\bigskip

\subsection{Proof of Theorem \ThConditionalA.iii when 
\poorBold{$\max_{0\leq i<n}g_{i/n}$} tends to infinity}
When $\max_{0\leq i<n}g_{i/n}$ tends to infinity, the proof of Theorem
\ThConditionalA\ requires an extra truncation, in part because the
function $v\in [0,\lambda)\mapsto (\lambda-v)^{\gamma-1}$ is no longer
bounded when $\gamma$ is less than $1$, and in part because the
various functionals introduced in the previous subsection are not well
behaved with respect to the convergence of measures.

To setup this truncation, we first prove the following lemma.

\Lemma{\label{iB}
  For any real number $B$ greater than $1$, define
  $$
    i_B(t)=\max\Bigl\{\, i\in \NN \,:\, {V(t)\over t} g_i > B \,\Bigr\} \, .
  $$

  \smallskip

  \noindent (i) If $\gamma$ is less than $1$, then 
  $$
    i_B(t)\sim (B/\gamma)^{1/(\gamma-1)} V(t)
    \qquad\hbox{and}\qquad 
    \g{i_B(t)}\sim (B/\gamma)^{\gamma/(\gamma-1)} t
  $$ 
  as $t$ tends to infinity.

  \smallskip

  \noindent (ii) If $\gamma$ is $1$, then 
  $$
    i_B(t)=o\bigl(V(t)\bigr) \qquad\hbox{and}\qquad 
    \g{i_B(t)}=o(t)
  $$
  as $t$ tends to infinity.
}

\bigskip

\Proof (i) Lemma \BS\ and Theorem 1.5.3 in Bingham, Goldie and Teugels (1989)
show that the sequence $(g_n)_{n\geq 0}$ is asymptotically equivalent to a 
nonincreasing sequence. Then, the asymptotic equivalence for $i_B$ follows 
from \gEquiv, and that for $\g{i_B}$ follows from Lemma \BS.

(ii) When $\gamma$ is $1$, the standard assumption ensures that the 
sequence $(g_n)_{n\geq 1}$ is asymptotically
equivalent to a nonincreasing sequence. For any positive real
number $c$ fixed, \gEquiv\ shows that $g_{\lfloor cV(t)\rfloor}\sim t/V(t)$,
and this forces $i_B(t)$ to be negligible compared to $V(t)$. The assertion
on $\g{i_B(t)}$ follows from Lemma \BS.\hfill\qed

\bigskip

Define the function
$$
  h(B)=\cases{2(B/\gamma)^{1/(\gamma-1)} & if $\gamma<1$, \cr
              1/B                        & if $\gamma=1$. \cr}
$$
Lemma \iB\ implies that for any fixed $B$, the function $i_B$ is ultimately
less than $h(B)V$. One sees that when $\gamma$ is $1$, the inequality 
$i_B<h(B)V$ holds ultimately for any positive function $h$. Our choice of
$h(B)=1/B$ in this case is entirely arbitrary and any positive function which
tends to $0$ at infinity could be used in what follows.

For any positive real number $B$, define
$$
  \calS_{B,t}(\lambda)
  ={1\over V(t)} \sum_{0\leq i<\lambda V(t)}
  \Bigl({V(t)\over t}g_i\wedge B\Bigr) X_{\lfloor \lambda V(t)\rfloor -i} \, .
$$
Our next lemma shows that given that $M$ exceeds a large level $t$ the
process $\St$ is well approximated by $\calS_{B,t}$ provided that $B$
is large enough.

\Lemma{\label{LemmaTruncationB}
  For any positive $\Lambda$ and $\epsilon$,
  $$
    \limB\limsupt \P\{\, |\St-\calS_{B,t}|_{[0,\Lambda]}>\epsilon
                         \mid M>t\,\}
    =0 \, .
  $$
}

\Proof In this proof, it is convenient to extend the sequence $(X_i)_{i\geq 1}$
to a sequence $(X_i)_{i\in \ZZ}$ of independent and identically distributed
random variables. Moreover, we define the centered random variables
$Z_i=|X_i|-\E |X_i|$, $i\in\ZZ$.

We rewrite $(\St-\calS_{B,t})(\lambda)$ as
$$
  {1\over V(t)} \sum_{0\leq i<\lambda V(t)} \Bigl( {V(t)\over t}g_i-B\Bigr)_+
  X_{\lfloor\lambda V(t)\rfloor -i} \, .
$$
Since $\bigl({V(t)\over t}g_i-B\bigr)_+$ vanishes for $i$ greater 
than $i_B(t)$, Lemma \iB\ and the discussion which follows shows that 
for any fixed $B$ greater than $1$, for any $t$ large enough and any 
positive $\lambda$,
$$\eqalignno{
  |(\St-\calS_{t,B})|(\lambda)
  &{}\leq {1\over V(t)} \sum_{0\leq i<h(B)V(t)} {V(t)\over t} g_i 
    |X_{\lfloor\lambda V(t)\rfloor -i}|
  &\equa{EqLemmaTruncationBA} \cr
  &{}={1\over t} \sum_{0\leq i<h(B)V(t)} g_i Z_{\lfloor\lambda V(t)\rfloor-i}
   + {\E |X_1|\over t} \g{h(B)V(t)} \, .\cr
  }
$$
Let $\epsilon$ be a positive real number. Since
$$
  {\E |X_1|\over t} \g{h(B)V(t)}\sim \E |X_1| h(B)^\gamma
$$
as $t$ tends to infinity, we take $B$ so large that 
$t^{-1}\E |X_1|\g{h(B)V(t)}$ is less than $\epsilon/2$ ultimately. 
Then \EqLemmaTruncationBA\ shows that there exists $t_0$, which does not
depend on $\lambda$, such that for any $t$ at least $t_0$ the probability
that $|\St-\calS_{t,B}|(\lambda)$ exceeds $2\epsilon$ is at most
$$
  \P\Bigl\{ \sum_{0\leq i<h(B)V(t)} g_iZ_{\lfloor\lambda V(t)\rfloor-i} >
  {t\over \g{h(B)V(t)}}\epsilon \g{h(B)V(t)}\,\Bigr\} \, .
  \eqno{\equa{EqLemmaTruncationBB}}
$$
Define $J_*$ as $J$ but substituting the distribution of $Z_i$ for that of 
\hbox{$X_i-\E X_i$}. Note that $t/\g{h(B)V(t)}\sim h(B)^{-\gamma}$. Then Proposition
\LDP\ implies that the logarithm of \EqLemmaTruncationBB\ is asymptotically
equivalent to $V(t)$ times
$$
  -h(B) J_*\Bigl({\epsilon\over h(B)^\gamma\gamma}\Bigr)
  \eqno{\equa{EqLemmaTruncationBC}}
$$
as $t$ tends to infinity. Since the logarithmic tail of the distribution 
function of $Z_i$ is regularly varying of index $\alpha$ greater than $1$,
so is $J_*$. Thus, given our choice of $h$, we can find $B$ large enough
so that the negative of \EqLemmaTruncationBC\ is greater than $3\theta$ 
---~recall that  $\theta$ was defined in section 2. Thus, for $t$ large 
enough, \EqLemmaTruncationBB\ is at most $\exp \bigl( -2\theta V(t)\bigr)$. 
Since 
$$
  |\St-\calS_{B,t}|_{[0,\Lambda]}
  = \max_{0\leq i<\Lambda V(t)} |\St-\calS_{B,t}|\Bigl({i\over V(t)}\Bigr)
  \, ,
$$
Bonferroni's inequality implies
$$
  \P\{\, |\St-\calS_{B,t}|_{[0,\Lambda]}>\epsilon\,\}
  =o(\P\{\, M>t\,\})
$$
as $t$ tends to infinity, and this proves the lemma.\hfill\qed

\bigskip

To prove Theorem \ThConditionalA\ under assumption (ii) of Theorem \ThCaseA, we
mostly repeat its proof under the assumption (i) of Theorem \ThCaseA, 
substituting $(\lambda-v)^{\gamma-1}\wedge B$ for $(\lambda -v)^{\gamma-1}$
and substituting the bound $(\lambda -v)^{\gamma-1}\wedge B\leq B$
for the bound $(\lambda-v)^{\gamma-1}\leq \Lambda^{\gamma-1}$. We indicate
the changes that are needed, from which it should be clear how the arguments
need to be changed.

Instead of the functional $\fS(\nu)$, we define
$$
  \fS_B(\nu)(\lambda)
  =\int\II_{[0,\lambda)}(v) \gamma\bigl((\lambda-v)^{\gamma-1}\wedge B\bigr)
  x\d\nu (x) \, ,
$$
and instead of the functional $\fS(\nu,b)$, we define
$$
  \fS_B(\nu,b)(\lambda)=\int\II_{[0,\lambda)}(v)
  \gamma\bigl((\lambda-v)^{\gamma-1}\wedge B\bigr) \sign(x)(|x|\wedge b)
  \d\nu(v,x) 
  \, .
$$
In what follows, we state the analogues of the lemmas of 
subsection \fixedref{6.3}. We do not indicate the proof when
it is identical to that of the previous subsection up to the substitutions
indicated above.

Our first lemmas are the analogues of Lemmas \LemmaApproxA, \LemmaApproxB,
and \LemmaApproxC.

\Lemma{
  For any positive $\Lambda$ and $\epsilon$,
  $$
    \limB \limsupt \P\bigl\{\, 
    |\calS_{B,t}-\fS_B(\calM_t)|_{[0,\Lambda]}>\epsilon 
    \bigm| M>t\,\bigr\}
    = 0 \, .
  $$%
}%
\finetune{\vskip -5pt}
\hfuzz=3pt
\Lemma{
  For any positive $\Lambda$ and $\epsilon$,
  $$\displaylines{
    \limB\limsup_{b\to\infty}\limsupt
    \P\bigl\{\, |\fS_B(\calM_t)-\fS_B(\calM_t,b)|_{[0,\Lambda]}>\epsilon 
    \bigm| M>t\,\bigr\} 
    \hfill\cr\hfill 
    = 0 \, .\cr}
  $$%
}%
\finetune{\vskip -15pt}
\hfuzz=0pt%
\Lemma{
  For any positive $\Lambda$,
  $$
    \limB\limsup_{b\to\infty} |\fS_B(\calM,b)-\fS_B(\calM)|_{[0,\Lambda]}
    =0 \, .
  $$
}

Instead of the functional $\fS_\epsilon(\nu,b)$, we define
$$
  \fS_{B,\epsilon}(\nu,b)(\lambda)
  = \int\Ieps(v/\lambda)\gamma\bigl((\lambda-v)^{\gamma-1}\wedge B\bigr)
    \sign(x)(|x|\wedge b) \d\nu(v,x) \, .
$$

\Lemma{
  For any positive $\Lambda$ and any $t$ large enough, both 
  $|\fS_{B,\epsilon}(\calM_t,b)-\fS_B(\calM_t,b)|_{[0,\Lambda]}$ and 
  $|\fS_{B,\epsilon}(\calM,b)-\fS_B(\calM,b)|_{[0,\Lambda]}$ are bounded
  by $2\gamma Bb\epsilon$.
  }

\bigskip

Referring to Lemma \LemmaApproxDA, instead of considering the modulus of
continuity of $\St$, we consider that of $\calS_{B,t}$. With an obvious
notation, we have the following.

\Lemma{
  For any positive $\epsilon$, $\Lambda$ and $B$,
  $$
    \lim_{\delta\to 0} \limsupt \P\{\, |\omega_{B,t,\delta}|_{[0,\Lambda]}
    >\epsilon \mid M>t\,\} = 0 \, .
  $$
}

Lemmas \LemmaApproxE\ and \LemmaApproxF\ and the conclusion of the proof do
not depend on $\gamma$. This proves Theorem \ThConditionalA\ under assumption
(ii) of Theorem \ThCaseA.

\bigskip

\subsection{Scaling argument}
We proved Theorem \ThConditionalA\ when $\mu$ is $-1$. To allow for
other negative values, as we did at the end of section \fixedref{5.2},
we index relevant quantities by the mean $\mu$ in parentheses so as to
make more transparent the scaling properties of various expressions.

Considering the innovations 
of the process, we first set $X_{(\mu),i}=(-\mu)X_{(-1),i}$, $i\geq 1$, so that
$M_{(\mu)}=(-\mu)M_{(-1)}$. Thus, the conditional probability given 
$M_{(\mu)}$ exceeds $t$ is the conditional probability given $M_{(-1)}$ 
exceeds $t/({-\mu})$. 

The random variables
$$
  \calN_{(\mu),t}
  = {V(t/{-\mu})\over V(t)} \calN_{(-1),t/{-\mu}}
$$
converges to $(-\mu)^{-1/\gamma}\tau_{(-1)}$ given $M_{(-1)}>t/{-\mu}$
as $t$ tends to infinity.
Referring to the how $\tau$ is defined in \HypTauUnique, equality \JScale\ 
implies
$$
  \tau_{(\mu)}=(-\mu)^{-1/\gamma}\tau_{(-1)} \, .
  \eqno{\equa{TauScale}}
$$
Therefore, as $t$ tends to infinity, $\calN_{(\mu),t}$ converges 
to $\tau_{(\mu)}$ given that $M_{(\mu)}$ exceeds $t$.

We write $\calM_{(\mu),t}$ as
$$\eqalign{
  \calM_{(\mu),t} 
  &{}={V(t/{-\mu})\over V(t)} {1\over V(t/{-\mu})} \sum_{i\geq 1} 
    \delta_{ ({V(t/{-\mu})\over V(t)} {i\over V(t/{-\mu})},-\mu X_{(-1),i}) } 
    \cr
  &{}={V(t/{-\mu})\over V(t)} \int \delta_{ ({V(t/{-\mu})\over V(t)} v,-\mu x)}
    \d\calM_{(-1),t/{-\mu}}(v,x) \, . \cr
  }
$$
Given $M_{(-1)}>t/(-\mu)$, we proved that the measures $\calM_{(-1),t/{-\mu}}$ 
converge to $\calM_{(-1)}$. It follows that given $M_{(\mu)}>t$, the measures
$\calM_{(\mu),t}$ converge to
$$
  \calM_{(\mu)}=(-\mu)^{-1/\gamma} \int\delta_{((-\mu)^{-1/\gamma}v,-\mu x)}
  \d\calM_{(-1)}(v,x) \, .
  \eqno{\equa{MNewDef}}
$$
Thus we need to check that this definition of $\calM_{(\mu)}$ coincides with
that in \MDef. With $\calM_{(\mu)}$ defined as in \MNewDef, we have, for any
bounded and continuous function $f$ on the right half-space,
$$\displaylines{\qquad
  \calM_{(\mu)}f
  =(-\mu)^{-1/\gamma} \int f\bigl( (-\mu)^{-1/\gamma}v,-\mu x\bigr) \d
    \calM_{(-1)} (v,x)
  \hfill\cr\qquad
  \phantom{\calM_{(\mu)}f}
  {}=(-\mu)^{-1/\gamma} \int f\bigl( (-\mu)^{-1/\gamma}v,-\mu x\bigr) 
  \hfill\cr\hfill 
    {}\times{ \exp\bigl( A_{(-1)}k_\gamma(v/\tau_{(-1)})
      x\bigr)
     \over
      \varphi_{(-1)} \bigl( A_{(-1)} k_\gamma(v/\tau_{(-1)})\bigr)
    }
   \d v\d F_{(-1)}(x) \, . 
  \cr}
$$
The change of variable $w=(-\mu)^{-1/\gamma}v$ and equality \TauScale\ yield
$$\displaylines{
  \calM_{(\mu)}f
  = \int \E f(w,X_{(\mu),1})
  \hfill\cr\hfill
  \times
  { \exp\bigl( A_{(-1)}k_\gamma(w/\tau_{(\mu)})
      X_{(\mu),1}/(-\mu)\bigr)
     \over
      \varphi_{(-1)} \bigl( A_{(-1)}k_\gamma (w/\tau_{(\mu)}) \bigr)
  } \d w \, .
  \quad
  \equa{MfEq}
  }
$$
Equality \PhiScale\ implies
$$
  m_{(-1)}(\lambda) = {1\over {-\mu}}m_{(\mu)} \Bigl({\lambda\over{-\mu}}\Bigr)
  \, . 
  \eqno{\equa{mScale}}
$$
Thus, $A_{(-1)}$ being defined in \ADefA, we have,
$$\eqalign{
  \tau_{(-1)}^{-\gamma}
  &{}= \int_0^1 k_\gamma(u) m_{(-1)}\bigl(A_{(-1)}k_\gamma(u)\bigr)\d u \cr
  &{}={1\over {-\mu}} \int_0^1 k_\gamma(u) 
    m_{(\mu)}\Bigl( {A_{(-1)}k_\gamma(u)\over {-\mu}} \Bigr) \d u \, .
  }
$$
It then follows from \TauScale\ that
$$
  \tau_{(\mu)}^{-\gamma}
  = \int_0^1 k_\gamma(u) m_{(\mu)}
  \Bigl( {A_{(-1)}k_\gamma(u)\over -\mu} \Bigr) \d u \, .
$$
Given \ADefA, this implies
$$
  A_{(\mu)}=A_{(-1)}/(-\mu) \, .
  \eqno{\equa{AScale}}
$$
Thus, referring to \MfEq\ and using \PhiScale, we obtain 
that $\calM_{(\mu)}f$ is equal to
$$
  \int \E f(w,X_{(\mu),1})
  { \exp\bigl( A_{(\mu)}k_\gamma(w/\tau_{(\mu)}) X_{(\mu),1} \bigr)
    \over
    \varphi_{(\mu)}\bigl( A_{(\mu)} k_\gamma(w/\tau_{(\mu)}) \bigr) 
  }
  \d w \, .
$$
This shows that definition \MNewDef\ for $\calM_{(\mu)}$ coincides with that
in \MDef.

We finally prove the convergence of the process $\calS_{(\mu),t}$. Since
$$\eqalign{
  \calS_{(\mu),t}(\lambda)
  &{}= { S_{(\mu),\lfloor\lambda V(t)\rfloor} \over t } \cr
  &{}= {-\mu S_{(-1),\lfloor\lambda V(t/-\mu) V(t)/V(t/-\mu)\rfloor} \over 
         (-\mu) t/(-\mu) } \, , \cr
  }
$$
the process $\calS_{(\mu),t}$ converges to 
$\calS_{(-1)}((-\mu)^{1/\gamma}\,\cdot\,)$ given $M_{(-1)}>t/{-\mu}$,
as $t$ tends to infinity.. Using the
definition of $\calS_{(-1)}$ in \SDef, the limiting process at $\lambda$ is
$$\displaylines{
 \calS_{(-1)}\bigl((-\mu)^{1/\gamma}\lambda\bigr)
  =\int_0^{\lambda (-\mu)^{1/\gamma}} 
  \gamma\bigl( \lambda (-\mu)^{1/\gamma}-v\big)^{\gamma-1} 
  \hfill\cr\hfill
  m_{(-1)}
  \bigl( A_{(-1)}k_\gamma(v/\tau_{(-1)})\bigr) \d v \, .
  \cr}
$$
The change of variable $v=(-\mu)^{1/\gamma}w$, equalities \mScale, \TauScale,
and \AScale\ yield that the limiting process at $\lambda$ is
$$
  \int_0^{\lambda} 
  \gamma(\lambda-w)^{\gamma-1} m_{(\mu)}
  \bigl( A_{(\mu)}k_\gamma (w/\tau_{\mu}) 
  \bigr) \d w \, ,
$$
which matches the definition of $\calS_{(\mu)}$ in \SDef. This proves Theorem
\ThConditionalA\ for arbitrary means.\hfill\qed

\bigskip


\noindent{\bf References}
\medskip

{\leftskip=\parindent \parindent=-\parindent
 \par

R.R.~Bahadur (1971). {\sl Some Limit Theorems in Statistics}, SIAM.

G.~Balkema, P.\ Embrechts (2007). {\sl High Risk Scenarios and Extremes, a
Geometric Approach}, European Mathematical Society.

Ph.~Barbe, M.~Broniatowski (1998). Note on functional large deviation principle
for fractional ARIMA processes, {\sl Statist.\ Inf.\ Stoch.\ Proc.}, 1, 17--27.

Ph.~Barbe, W.P.~McCormick (2008). Veraverbeke's theorem at large ---
on the maximum of some processes with negative drifts and heavy tail
innovations, preprint ({\tt arXiv:0802.3638}).

O.E.~Barndorff-Nielsen (1978). {\sl Information and Exponential Families in
Statistical Theory}, Wiley.

P.~Billingsley (1968). {\it Convergence of Probability Measures}, Wiley.

N.H.~Bingham, C.M.~Goldie, J.L.~Teugels (1989). {\sl Regular Variation}, 
2nd ed., Cambridge University Press.

N.H.~Bingham, J.~Teugels (1975). Duality for regularly varying functions,
{\sl Quarterly J.\ Math.}, 26, 333--353.

M.~Broniatowski, A.~Fuchs (1995). Tauberian theorems, Chernoff inequality,
and the tail behavior of finite convolutions of distribution functions,
{\sl Adv.\ Math.}, 116, 12-33.

L.D.~Brown (1986). {\sl Fundamentals of Statistical Exponential Families with
Applications in Statistical Decision Theory}, IMS.

R.~Burton, H.~Dehling (1990). Large deviations for some weakly dependent
random processes, {\sl Statist.\ Probab.\ Lett.}, 9, 397--401.

C.-S.~Chang, D.D.~Yao, T.~Zajic (1999). Large deviations, moderate 
deviations, and queues with long-range dependent input, {\sl Adv.\ Appl.\ 
Probab.}, 31, 254--277.

J.~Collamore (1996). Hitting probabilities and large deviations, {\sl Ann.\
Probab.}, 24, 2065--2078.

J.~Collamore (1998). First passage times of general sequences of random
vectors: a large deviations approach, {\sl Stoch.\ Proc.\ Appl.}, 78, 97--130.

I.~Csisz\'ar (1984). Sanov property, generalized $I$-projection and a 
conditional limit theorem, {\sl Ann.\ Probab.}, 12, 768--793.

A.~Dembo, O.~Zeitouni (1993). {\sl Large Deviations Techniques and 
Applications}, Jones and Bartlett.

P.~Diaconis, D.~Freedman (1988). Conditional limit theorems for exponential
families and finite version of de Finetti's theorem, 
{\sl J.\ Theoret.\ Probab.}, 1, 381--410.

A.B.~Dieker, M.~Mandjes (2006). Efficient simulation of random walks exceeding
a nonlinear boundary, {\sl Stoch.\ Models}, 22, 459--481.

N.G.~Duffield, N.~O'Connell (1995). Large deviations and overflow probabilities
for the general single-server queue, with applications, {\sl Math.\ Proc.\
Camb.\ Phil.\ Soc.}, 118, 363--374.

N.G.~Duffield, W.~Whitt (1998). Large deviations of inverse processes with
nonlinear scalings, {\sl Ann.\ Appl.\ Probab.}, 4, 995--1026.

R.~Ellis (1984). Large deviations for a general class of random vectors,
{\sl Ann.\ Probab.}, 12, 1--12.

W.~Feller (1971). {\sl An Introduction to Probability Theory and its 
Applications}, Wiley.

M.I.~Freidlin, A.D.~Wentzell (1984). {\sl Random Perturbation of Dynamical
Systems}, Springer.

A.J.~Ganesh, N.~O'Connell (2002). A large deviation principle with queueing
applications, {\sl Stochastics and Stochastic Reports}, 73, 25--35.

J.~G\"artner (1977). On large deviations from invariant measure, {\sl Theor.\
Probab.\ Appl.}, 22, 24--39.

H.U.~Gerber (1982). Ruin theory in the linear model, {\sl Insurance: 
Mathematics and Econonmics}, 1, 177-184.

P.~Groeneboom, J.~Oosterhoff, F.~Ruymgaart (1979). Large deviation
theorems for empirical probability measures, {\sl Ann.\ Probab.}, 7, 553--586.

J.M.~Hammersley, D.G.~Handscomb (1964). {\sl Monte Carlo Methods}, Chapman \& 
Hall.

A.B.~Hoadley (1967). On the probability of large deviations of functions of
several empirical cdf's, {\sl Ann.\ Math.\ Statist.}, 38, 360--381.

J.~H\"usler, V.~Piterbarg (2004). On the ruin probability for physical 
fractional Brownian motion, {\sl Stoch.\ Proc.\ Appl.}, 113, 315--332.

J.\ Iscoe, P.\ Ney, E.\ Nummelin (1985). Large deviations of uniformly
recurrent Markiv additive processes, {\sl Adv.\ Appl.\ Math.}, 6, 373-412.

J.~Janssen (1982). On the interaction between risk and queueing theories,
{\sl Bl\"atter der DGVFM}, 15, 383--395.

Y.~Kasahara (1978). Tauberian theorems of exponential type, {\sl J.\ Math.\
Kyoto Univ.}, 18, 209--219.

H.~Kesten (1973). Random difference equations and renewal theory for product
of random matrices, {\sl Acta. Math.}, 131, 207--248.

G.~Letac (1992). {\sl Lectures on Exponential Families and their Variance
Functions}, IMPA.

T.~Lindvall (1973). Weak convergence of probability measures and random
functions in the function space $\D[\,0,\infty)$, {\sl J.\ Appl.\ Probab.},
10, 109--121.

A.A.~Mogulskii (1976). Large deviations for trajectories of multi-dimensional
random walks, {\sl Theor.\ Probab.\ Appl.}, 21, 300--315.

A.~M\"uller, G.~Pflug (2001). Asymptotic ruin probabilities for risk processes
with dependent increments, {\sl Insurance Math.\ Econom.}, 28, 381--392.

H.~Nyrhinen (1994). Rough limit results for level-crossing probabilities,
{\sl J.\ Appl.\ Probab.}, 31, 373--382.

H.~Nyrhinen (1995). On the typical level crossing time and path, {\sl Stoch.\
Proc.\ Appl.}, 58, 11--137.

H.~Nyrhinen (1998). Rough descriptions of ruin for a general class of surplus
processes, {\sl Adv.\ Appl.\ Probab.}, 30, 1008--1026.

N.U.~Prabhu (1961). On the ruin problem of collective risk theory, {\sl Ann.\
Math.\ Statist.}, 32, 757--764.

S.D.~Promislow (1991). The probability of ruin in a process with dependent
increments, {\sl Insurance Math.\ Econom.}, 10, 99--107.

R.T.~Rockafellar (1970). {\sl Convex Analysis}, Princeton University Press.

W.~Rudin (1976). {\sl Principle of Mathematical Analysis}, 3rd ed., 
McGraw-Hill.

J.~Sadowsky (1996). On Monte Carlo estimation of large deviations 
probabilities, {\sl Ann.\ Appl.\ Probab.}, 6, 399--422.

}


\bigskip

\setbox1=\vbox{\halign{#\hfil&\hskip 40pt #\hfill\cr
  Ph.\ Barbe            & W.P.\ McCormick\cr
  90 rue de Vaugirard   & Dept.\ of Statistics \cr
  75006 PARIS           & University of Georgia \cr
  FRANCE                & Athens, GA 30602 \cr
                        & USA \cr
                        & bill@stat.uga.edu \cr}}
\box1

\bye